\documentclass[11pt]{amsart}   	
\usepackage{geometry}             
\usepackage[dvipsnames]{xcolor}   		             		
\usepackage{graphicx}				
\usepackage{amssymb,mathrsfs}
\usepackage{amssymb,amsthm,amsmath,stmaryrd}
\usepackage{tikz}
\usepackage{tikz-cd}
\usepackage{accents,upgreek,enumerate}
\usepackage[headings]{fullpage}
\usepackage{bm}
\usepackage[all]{xy}
\usepackage{caption}
\usepackage{csquotes}
\usepackage{enumitem}
\usepackage{comment}

\usetikzlibrary{calc}
\usetikzlibrary{fadings}
\usetikzlibrary{decorations.pathmorphing}
\usetikzlibrary{decorations.pathreplacing}

\renewcommand{\familydefault}{ppl}
\DeclareMathAlphabet{\mathpzc}{OT1}{pzc}{m}{it}

\usepackage[backref]{hyperref}
\hypersetup{
colorlinks   = true,          %Colors links instead of ugly boxes
urlcolor     = violet,          %Color for external hyperlinks
linkcolor    = Purple,          %Color of internal links
citecolor   = RubineRed             %Color of citations
}

\newtheorem{theorem}{Theorem}[section]
\newtheorem{corollary}[theorem]{Corollary}
\newtheorem{lemma}[theorem]{Lemma}
\newtheorem{proposition}[theorem]{Proposition}

\newtheorem{quasi-theorem}[theorem]{Quasi-Theorem}

\theoremstyle{definition}
\newtheorem{definition}[theorem]{Definition}

\newtheorem{example}[theorem]{Example}
\newtheorem{remark}[theorem]{Remark}

\newcommand{\CC} {{\mathbb C}}

\newcommand{\RR} {{\mathbb R}}		
\newcommand{\ZZ} {{\mathbb Z}}	
\renewcommand{\AA}{\mathbb{A}}	

\def\setminus{\smallsetminus}

\newcommand{\Puiseux}[2][t]{#2\{\!\{#1\}\!\}}

\newcommand{\Tr}{\operatorname{tr}}

\DeclareMathOperator{\chara}{char} 
 
\DeclareMathOperator{\conv}{conv}

\DeclareMathOperator{\Area}{Area}

\DeclareMathOperator{\val}{val}
\DeclareMathOperator{\floor}{floor}

\DeclareMathOperator{\GW}{GW}

\newcommand{\bw}{\mathbf{w}}

\newcommand{\cal}{\mathcal}

\def\cD{{\cal D}}
\def\cE{{\cal E}}

\def\trop{\mathrm{trop}}
\def\floor{\mathrm{floor}}
\def\mult{\mathrm{mult}}

\makeatletter
\def\blfootnote{\xdef\@thefnmark{}\@footnotetext}
\makeatother

\title{\bf Arithmetic counts of tropical plane curves and their properties} 

\date{}

\author {Andr\'es Jaramillo Puentes}
\address {Universit\"at Duisburg-Essen, Fakult\"at f\"ur Mathematik, Thea-Leymann-Str. 9, 45127 Essen, Germany}
\email {andres.jaramillo-puentes@uni-due.de}
\author {Hannah Markwig}
\address {Universit\"at T\"ubingen, Fachbereich Mathematik, Auf der Morgenstelle 10, 72076 T\"ubingen, Germany }
\email {hannah@math.uni-tuebingen.de}
\author {Sabrina Pauli}
\address {TU Darmstadt, Fachbereich Mathematik, Schlossgartenstraße 7, 64289 Darmstadt, Germany}
\email {pauli@mathematik.tu-darmstadt.de}
\author{Felix R\"ohrle}
\address {Universit\"at T\"ubingen, Fachbereich Mathematik, Auf der Morgenstelle 10, 72076 T\"ubingen, Germany }
\email{roehrle@math.uni-tuebingen.de}


\begin{document}

\begin{abstract}
%Tropical plane curves can be used to count rational algebraic plane curves satisfying point conditions both over the complex and over the real numbers. 
%In both cases, we consider the same set of tropical plane curves, depending on the degree and the choice of point conditions. We only vary the multiplicity with which each tropical curve is counted. The arithmetic or quadratically enriched count of plane curves can be specialized to arbitrary fields. 
Recently, the first and third author proved a correspondence theorem which recovers the Levine-Welschinger invariants of toric del Pezzo surfaces as a count of tropical curves weighted with arithmetic multiplicities. 
%When specializing to the complex or real numbers, we obtain the well-known tropical complex resp.\ real multiplicities. 
In this paper, we study properties of the arithmetic count of plane tropical curves satisfying point conditions.
We prove that this count is independent of the configuration of point conditions. Moreover, a Caporaso-Harris formula for the arithmetic count of plane tropical curves is obtained by moving one point to the very left. Repeating this process until all point conditions are stretched, we obtain an enriched count of floor diagrams which coincides with the tropical count. Finally, we prove polynomiality properties for the arithmetic counts using floor diagrams. 
%Via the correspondence theorem these results give new insight into the structure of enriched counts in algebraic geometry.
%To sum up, the arithmetic enrichment of the tropical plane curve count interacts nicely with its combinatorics.
\end{abstract}

\maketitle

\vspace{-0.2in}

\setcounter{tocdepth}{1}
\tableofcontents

\section{Introduction}

\newcommand{\counttheorems}{\thetheorem}
\renewcommand*{\thetheorem}{\Alph{theorem}}

Enumerative geometry poses geometric counting problems defined over some fixed base field $K$. Often, $K$ is taken to be the field of complex numbers since this field is algebraically closed, thus leading to desirable properties of the counting problem. For example, consider the count of complex plane curves of fixed degree $d$ and genus $g$ passing through some points. This number is a so-called Gromov-Witten invariant of the plane and is denoted $N(d,g)$. It is indeed an invariant in the sense that $N(d,g)$ does not depend on the position of the points, as long as they are chosen generically. In contrast, when we choose $K = \RR$ to count real plane curves satisfying real point conditions, the invariance does no longer hold. For example, there can be $8$, $10$, or $12$ rational real cubics passing through $8$ real points in the plane \cite{DK00}. However, over $\mathbb{R}$ one can often introduce signed counts which are invariant, and where the sign naturally depends on geometric properties of the objects to be counted. In the example of rational (i.e. genus $0$) plane curves subject to point conditions, Welschinger established an invariant $W(d)$ by counting each curve with a sign taking the number of solitary nodes into account \cite{Wel05}.  In higher genus, no analogue of a Welschinger invariant is known yet. 

Since recently, $\mathbb{A}^1$-enumerative geometry offers a universal approach to geometric counting problems by considering enumerative invariants valued in the Grothendieck-Witt ring $\GW(K)$ of the base field $K$. 
For example, there are arithmetic counts of lines in a cubic surface \cite{KW17}, bitangents to a smooth plane quartic curve \cite{LV19}, and plane curves satisfying point conditions \cite{LevineWelschinger,KassLevineSolomonWickelgren}.
In each of these examples, the complex and real invariants can be recovered from the $\AA^1$-count.

The success of tropical methods in enumerative geometry started with Mikhalkin's correspondence theorem~\cite{Mi03}, in which he proved that both  $N(d,g)$ and $W(d)$ can be computed as weighted counts of rational plane tropical curves of fixed degree satisfying point conditions in $\RR^2$. Remarkably,  in genus $0$, the two tropical counts are in fact summing over the same set of tropical curves, but counting them with complex or real multiplicity to obtain the Gromov-Witten or Welschinger invariant, respectively.
It is interesting to observe that the tropical count with real multiplicity readily generalizes to curves of higher genus, even though there is no analogue of Welschinger's numbers in higher genus. Despite these tropical invariants being in some sense artificial invariants on the tropical side, they have been used to obtain interesting recursions for the algebraic numbers \cite{IKS09}. Conversely, one can use the technique of the correspondence theorem to construct real algebraic curves tropicalizing to the real tropical curve count also in higher genus --- we call them: algebraic curves \emph{close to the tropical limit} --- but since there is no invariant on the real side, the gain of this direct correspondence is limited.

The tropical approach is in fact not restricted to the case of counting curves in the projective plane $\mathbb{P}^2$, instead we can count curves in any toric surface. A lattice polygon $\Delta$ defines a toric surface as well as a curve class given by hyperplane sections. We denote the count of complex curves of class $\Delta$ and genus $g$ satisfying point conditions by $N(\Delta,g)$. In the case of $\mathbb{P}^2$ and curves of degree $d$ (corresponding to the triangle $\Delta_d$ with vertices $(0,0)$, $(0,d)$ and $(d,0)$), we continue to use the simplified notation $N(d,g)$.

The tropical approach taken to arithmetic counts of bitangents \cite{MPS22} as well as to the B\'ezout's theorem, or more generally to the Bernstein–Kushnirenko theorem \cite{puentes2022quadratically}, reveals how well tropical degenerations interact with arithmetic multiplicities. More recently, the first and third author proved an enriched correspondence theorem \cite{JPP23} expressing Levine's enriched Welschinger count for rational curves on toric del Pezzo surfaces $N_{\AA^1}(\Delta, 0)$ \cite{LevineWelschinger} as a tropical arithmetic count $N_{\AA^1}^\trop (\Delta, 0)$. The latter again follows the well-known pattern of summing over tropical curves but this time using suitable arithmetic multiplicities. 

 In this paper, we introduce an (\enquote{artificial}) tropical invariant $N_{\AA^1}^\trop (\Delta, g)$ of plane tropical curves of any genus $g$ counted with arithmetic multiplicity. In the case of genus $0$ on a toric del Pezzo surface, the correspondence theorem from \cite{JPP23} states the equality of these numbers with their arithmetic counterparts.  In higher genus, our tropical invariants correspond to an arithmetic count of algebraic curves close to the tropical limit.

The contribution of this paper is the study of combinatorial properties of the tropical numbers $N_{\AA^1}^\trop (\Delta, g)$.   In genus $0$, in the del Pezzo case and if all unbounded ends of a tropical curve $C$ have weight 1, then the correspondence theorem~\cite{JPP23} dictates that the arithmetic multiplicity of a tropical curve $C$ is essentially a linear combination of the real and complex multiplicities of $C$. This readily implies that various known results on the structure of real and complex tropical curve counts carry over to these new arithmetic counts.  The same holds for the more general tropical invariants, i.e.\ for higher genus and other toric surfaces, as the multiplicity is defined analogously.
Nevertheless, we make it the purpose of this paper to reprove these results by repeating the key arguments from the existing literature and manually adapting them to arithmetic counting. We believe that this showcases that the tropical approach to enumerative geometry is in fact sufficiently robust to accommodate this vast generalization of curve counting, while also providing a comprehensive overview over combinatorial methods which have become standard in complex curve counting and which are now making their way into arithmetic enumerative geometry.

\subsection{Results}

Let $C$ be a tropical plane curve with Newton polygon $\Delta$. 
 Let $K$ be a field of large enough characteristic, i.e.\ larger than the diameter of $\Delta$.
Assume that the unbounded ends of $C$ have weights $w_1, \ldots, w_k$. In Definition~\ref{def-troparith2} we give a straight forward generalization of the arithmetic multiplicity $\mult_{\AA^1}(C)$ defined in \cite{JPP23} by incorporating the $w_i$. 
This leads us to introduce a tropical arithmetic count $N^{\trop}_{\mathbb{A}^1}(\Delta,g)$ for tropical curves of genus $g$ and degree $\Delta$. Our first result is the invariance, which we prove by combining the tools for showing the complex and the real invariance:

\begin{theorem}[Invariance of arithmetic counts of plane tropical curves, Theorem~\ref{thm-invariance}]
	 Let $K$ be a field of large enough characteristic, i.e.\ larger than the diameter of $\Delta$.
 The count $N^\trop_{\mathbb{A}^1}(\Delta,g)\in \GW(K)$ does not depend on the location of the points in $\RR^2$ (as long as they are in general position) through which we require the tropical curves to pass.
\end{theorem}

In the special case of genus $0$ curves on toric del Pezzo surfaces, this invariance follows \textit{a posteriori} from the enriched correspondence theorem \cite{JPP23}. We establish the invariance of the tropical count in full generality  (i.e.\ also for higher genus, arbitrary toric surfaces and any weights of ends) by combinatorial considerations. More precisely, the space of point configurations in $\RR^2$ has a wall-and-chamber structure and we show that $N^\trop_{\mathbb{A}^1}(\Delta,g)$ remains the same when crossing any walls.

Since the count is invariant, we can move one of our points very far to the left, leading to a further degeneration of tropical plane curves. To use this degeneration in a recursion, we first need to further enlarge our counting problem by allowing ends of higher weight for our tropical curves, corresponding to higher contact orders with the toric boundary. 
We spell this out for the case of $\mathbb{P}^2$, obtaining a recursive formula for $N^{\trop, \alpha, \beta}_{\mathbb{A}^1}(d,g)$, similar to the Caporaso-Harris formula for Gromov-Witten invariants \cite{CH98}. Our proof builds on the tropical proof of the Caporaso-Harris formula \cite{GM052}.

\begin{theorem}[Caporaso-Harris formula, details see Theorem~\ref{thm-caporaso-harris}]
 Let $K$ be a field of characteristic larger than $d$.
	Fix a degree $d > 1$, genus $g$, and numbers $\alpha$ of fixed ends and $\beta$ non-fixed ends. Then the tropical arithmetic count $ N_{\mathbb{A}^1}^{\trop,\alpha,\beta}(d,g)\in \GW(K)$  satisfies a Caporaso-Harris-style recursion formula.
\end{theorem}

Applying the step in the Caporaso-Harris recursion iteratively until all points are stretched in a small horizontal strip, the tropical curves to be counted can be controlled by even simpler combinatorial objects called \emph{floor diagrams} \cite{BM08, FM09}. In Definition~\ref{def-arithmultfloor} we introduce the arithmetic multiplicity $\mult_{\AA^1}(\mathcal{D})$ for a floor diagram $\mathcal{D}$. In order to keep the notation reasonably simple, we give our definition only for curves on Hirzebruch surfaces, i.e. with Newton polygon of the form $\Delta = \Delta_k(a,b)$, see Figure \ref{fig-Hirz}.

\begin{theorem}[Theorem~\ref{thm-floor=trop}]
	Fix a genus $g$ and degree for tropical curves with Newton polygon of the form $\Delta = \Delta_k(a,b)$. 
  Let $K$ be a field of large enough characteristic, i.e.\ larger than the diameter of $\Delta$. Then the arithmetic count of floor diagrams equals the arithmetic count of tropical curves,
	$N_{\mathbb{A}^1}^{\floor} \big(\Delta,g \big) = N_{\mathbb{A}^1}^{\trop} \big(\Delta,g \big)  \in \GW(K).$
\end{theorem}

Finally, we use the floor diagram count of arithmetic numbers to prove polynomiality properties. 
The investigation of such properties goes back to the case of double Hurwitz numbers \cite{gjv:ttgodhn}.
% which count genus $g$ degree $d$ covers of $\mathbb{P}^1$ with two special ramification profiles $\mu,\nu \vdash d$ over $0$ and $\infty$ and only simple ramification else, at fixed points. Goulden, Jackson, and Vakil showed that these numbers are polynomial in the entries $\mu_i,\nu_j$ of the two partitions $\mu,\nu$ , and related this to interesting geometric properties of counts of covers. 
The piecewise polynomial structure of double Hurwitz numbers can easily be obeserved tropically \cite{cjm:wcfdhn}. The similarity of the count of tropical covers and the count of floor diagrams counting tropical plane curves made it possible to apply these tropical methods to obtain a piecewise polynomiality result for counts of curves in Hirzebruch surfaces \cite{AB14}. 
More precisely, we fix the width $a$ of the polygon in Figure \ref{fig-Hirz}, and we fix the number of left and right ends respectively. We consider the tuple of weights $\bw^{(r)}$ of the right ends and $\bw^{(l)}$ of the left ends as variables. To avoid case distinction we assume all weight to be odd.
%, where $\nu$ is a partition of $b$ and $\mu$ of  $ ak+b$. 
Then the complex count of tropical curves (and via the correspondence theorem, the complex count) can be viewed as a function in the variables $(\bw^{(r)},\bw^{(l)})$. 
%As the arithmetic multiplicity can be viewed as an interpolation between the complex and real count, one can deduce also piecewise polynomiality properties of the arithmetic count.
Using polynomiality of the complex count we show:

\begin{theorem}[Theorem~\ref{thm-piecepoly}]
	Let $K$ be a field of characteristic $0$.
    The function $N_{\mathbb{A}^1}^{\floor}(k,a,g)$ defined in Equation~\eqref{eq:floor_diag_count} associates to a tuple $(\bw^{(r)},\bw^{(l)})$ of odd weights an element in $\GW(K)$ which equals a multiple of $\mathbb{H}$ plus a sum of quadratic forms of the form $\big\langle (-1)^i\cdot \prod\bw^{(r)}_i\prod\bw^{(l)}_j \big\rangle$. The coefficient of  $\mathbb{H}$ is piecewise polynomial in the entries of $\bw^{(r)}$ and $\bw^{(l)}$.
\end{theorem}

Finally, fixing the cogenus $\delta$ (i.e. the number of nodes) and keeping the degree $d$ as a variable, \cite{FM09} have shown that the complex count is eventually polynomial. We adapt their proof to arithmetic counts.

\begin{theorem}[Theorem~\ref{thm-polynomiality}]
    Let $K$ be a field of characteristic $0$.
	For every fixed number of nodes $\delta$ there are polynomials $P, Q \in \mathbb{Q}[d]$ of degree $2\delta$ such that for all sufficiently large $d$ we have 
	\[ N_{\mathbb{A}^1}^\delta(d) = P(d) \mathbb{H} + Q(d) \langle 1 \rangle. \]
\end{theorem}

\subsection{Acknowledgements}
We thank Johannes Rau for useful discussions. We thank two anonymous referees for useful comments on how to improve the exposition.
The second and fourth author acknowledge support by DFG-grant MA 4797/9-1. 
The first and third author have been supported by the ERC programme QUADAG.  This paper is part of a project that has received funding from the European Research Council (ERC) under the European Union's Horizon 2020 research and innovation programme (grant agreement No. 832833).\\ 
\includegraphics[scale=0.08]{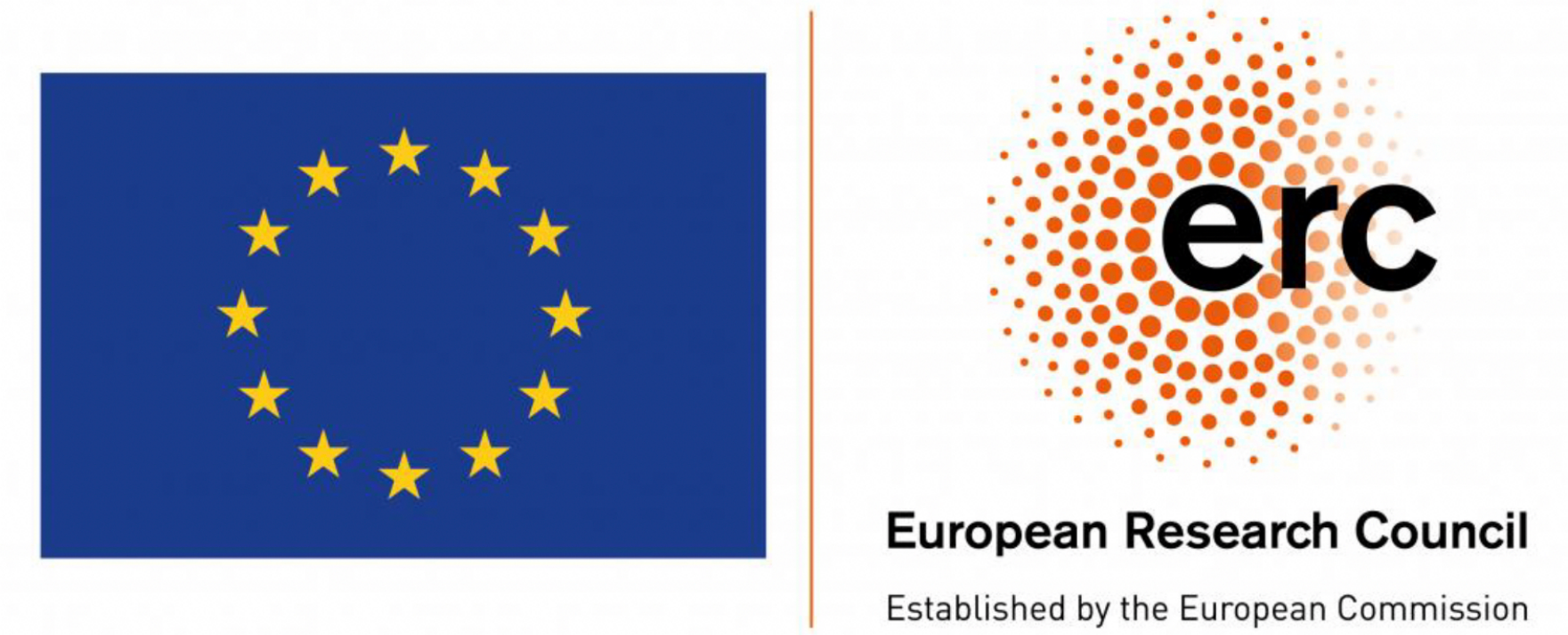}

\renewcommand{\thetheorem}{\arabic{section}.\arabic{theorem}}

\section{Preliminaries}
\subsection{Tropical plane curves}

We start with a brief overview of the theory of tropical plane curves. Readers familiar with this subject may jump directly to Section \ref{subsec-counts}. The exposition in this subsection is based on Chapter 3 in \cite{CMR23}. More details can be found there.

Consider a field $K$. A plane algebraic curve (or, more generally, a curve embedded in a toric surface) over $K$ is given as the zero-set of a polynomial in $K[x, y]$. In tropical geometry we can use the analogous definition, by just adding the word \enquote{tropical} everywhere:
tropical plane curves are tropical vanishing loci of tropical polynomials over the tropical semifield. As this perspective is not of particular relevance for our purpose, we just refer interested readers to surveys on tropical geometry such as \cite{RGST05, BS14}. 

Alternatively, one can also start with a field $K$ which is equipped with a non-Archimedean valuation $\val : K^\times \to \RR$ and an algebraic curve $\mathcal{C}$ defined over $K$. A tropical plane curve is the closure of the image of $\mathcal{C} \cap (K^\times)^2$ under the \emph{tropicalization map}
\[ \trop : (K^\times)^2 \longrightarrow \RR^2, \quad (x,y) \longmapsto \big( -\val(x), -\val(y) \big) \ . \]
The resulting tropical curve is called the \emph{tropicalization} of $\mathcal{C}$. It is a piecewise linear graph with edges of rational slope in $\mathbb{R}^2$. In fact, the two definitions mentioned here are equivalent by Kapranov's theorem, see \cite{RGST05, BS14}.

Every tropical plane curve $C$ is dual to a Newton subdivision of the Newton polygon of its defining equation defined as follows. Let $C$ be given by either a tropical polynomial $F$ or by a polynomial $f$ over the non-Archimedean field $(K, \val)$. Consider the Newton polygon of $F$ resp.\ $f$, and extend it into 3-space by adding to every lattice point the tropical coefficient of the corresponding monomial in $F$, resp.\ the negative of the valuation of the coefficient of $f$ as third coordinate. We can project the lower convex hull of this extended Newton polytope from above and obtain the dual Newton subdivision. Here, duality is in the sense that vertices of $C$ correspond to polygons of the subdivision, edges of $C$ to edges of the subdivision in orthogonal direction, and connected components of $\mathbb{R}^2\setminus C$ to vertices of the subdivision.
For more details on this construction, we refer again to \cite{RGST05, BS14}.

\begin{example}\label{ex-dual}
Let $K$ be the field of Puiseux series $\Puiseux{\CC}$ whose valuation sends a Puiseux series to its least exponent, i.e.\ its order.
Consider the degree $2$ polynomial $f \in \Puiseux{\CC}[x, y]$ given as
$$f=1+ t^{-1}\cdot  x+ x^2+  t^{-1}\cdot y+  t^{-1}\cdot x y+ y^2.$$
The tropical plane curve defined by $f$ is depicted in Figure \ref{fig:dualcurve}.
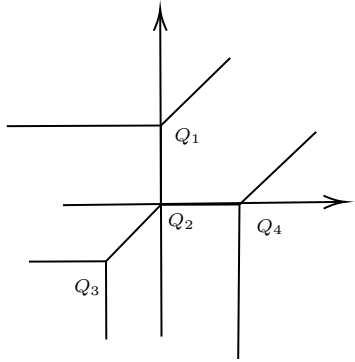
\begin{figure}
	
	\begin{center}
		\tikzset{every picture/.style={line width=0.75pt}} %set default line width to 0.75pt        

\begin{tikzpicture}[x=0.75pt,y=0.75pt,yscale=-1,xscale=1]
	%uncomment if require: \path (0,784); %set diagram left start at 0, and has height of 784
	
	%Straight Lines [id:da6208660064350149] 
	\draw    (122,328.8) -- (149.7,300.27) ;
	%Straight Lines [id:da776313181604554] 
	\draw    (149.7,260.27) -- (149.7,300.27) ;
	%Straight Lines [id:da21256404137119744] 
	\draw    (149.7,300.27) -- (189.33,300.13) ;
	%Straight Lines [id:da09907149260319081] 
	\draw    (189.33,300.13) -- (228,263.47) ;
	%Straight Lines [id:da3402388277882751] 
	\draw    (189.33,300.13) -- (188.68,378.1) ;
	%Straight Lines [id:da041298461105825135] 
	\draw    (71.95,260.6) -- (149.7,260.27) ;
	%Straight Lines [id:da3904803443612095] 
	\draw    (83.12,328.96) -- (122,328.8) ;
	%Straight Lines [id:da09098394786160224] 
	\draw    (122,328.8) -- (122.16,368.27) ;
	%Straight Lines [id:da5317229073223257] 
	\draw    (184.67,226.13) -- (149.7,260.27) ;
	%Straight Lines [id:da46367867974592203] 
	\draw    (150,366.8) -- (149.34,203.47) ;
	\draw [shift={(149.33,201.47)}, rotate = 449.77] [color={rgb, 255:red, 0; green, 0; blue, 0 }  ][line width=0.75]    (10.93,-3.29) .. controls (6.95,-1.4) and (3.31,-0.3) .. (0,0) .. controls (3.31,0.3) and (6.95,1.4) .. (10.93,3.29)   ;
	%Straight Lines [id:da29709890655713234] 
	\draw    (100.33,300.47) -- (240.67,298.82) ;
	\draw [shift={(242.67,298.8)}, rotate = 539.3299999999999] [color={rgb, 255:red, 0; green, 0; blue, 0 }  ][line width=0.75]    (10.93,-3.29) .. controls (6.95,-1.4) and (3.31,-0.3) .. (0,0) .. controls (3.31,0.3) and (6.95,1.4) .. (10.93,3.29)   ;
	
	% Text Node
	\draw (155,260.13) node [anchor=north west][inner sep=0.75pt]  [font=\tiny] [align=left] {$\displaystyle Q_{1}$};
	% Text Node
	\draw (151.7,303.27) node [anchor=north west][inner sep=0.75pt]  [font=\tiny] [align=left] {$\displaystyle Q_{2}$};
	% Text Node
	\draw (104.37,336.6) node [anchor=north west][inner sep=0.75pt]  [font=\tiny] [align=left] {$\displaystyle Q_{3}$};
	% Text Node
	\draw (196.37,305.94) node [anchor=north west][inner sep=0.75pt]  [font=\tiny] [align=left] {$\displaystyle Q_{4}$};

\end{tikzpicture}
	\end{center}

	 \caption{The tropical curve given by the polynomial in Example \ref{ex-dual}. It is dual to the subdivision depicted on the right of Figure \ref{fig:subdivision}.}
    \label{fig:dualcurve}
\end{figure}

The  pictures in Figure \ref{fig:subdivision} show the Newton polygon of $f$, i.e.\ the convex hull of the exponent vectors. Next to each point, we also note the height (i.e.\ the negative of the valuation of the coefficient) which is used to produce the extended Newton polytope in $3$-space. We project the upper faces down to the Newton polygon. We can see that the points of height $1$ form a triangular upper face of the extended Newton polytope. Two points of height $1$ together with a neighbour of height $0$ form a triangular face which is bent slightly downwards from this top part of the roof. Projecting to the Newton polygon, we obtain the picture on the right of Figure \ref{fig:subdivision}.

\begin{figure}
	\begin{center}
	\tikzset{every picture/.style={line width=0.75pt}} %set default line width to 0.75pt        

\begin{tikzpicture}[x=0.75pt,y=0.75pt,yscale=-1,xscale=1]
	%uncomment if require: \path (0,784); %set diagram left start at 0, and has height of 784
	
	%Straight Lines [id:da6208660064350149] 
	\draw    (102.35,38.97) -- (179.46,116.61) ;
	%Shape: Ellipse [id:dp45604016944512815] 
	\draw  [fill={rgb, 255:red, 0; green, 0; blue, 0 }  ,fill opacity=1 ] (98.17,38.97) .. controls (98.17,36.65) and (100.04,34.78) .. (102.35,34.78) .. controls (104.66,34.78) and (106.53,36.65) .. (106.53,38.97) .. controls (106.53,41.28) and (104.66,43.16) .. (102.35,43.16) .. controls (100.04,43.16) and (98.17,41.28) .. (98.17,38.97) -- cycle ;
	%Shape: Ellipse [id:dp25566478646064195] 
	\draw  [fill={rgb, 255:red, 0; green, 0; blue, 0 }  ,fill opacity=1 ] (97.85,78.12) .. controls (97.85,75.8) and (99.72,73.92) .. (102.03,73.92) .. controls (104.33,73.92) and (106.21,75.8) .. (106.21,78.12) .. controls (106.21,80.43) and (104.33,82.31) .. (102.03,82.31) .. controls (99.72,82.31) and (97.85,80.43) .. (97.85,78.12) -- cycle ;
	%Shape: Ellipse [id:dp5755207745665273] 
	\draw  [fill={rgb, 255:red, 0; green, 0; blue, 0 }  ,fill opacity=1 ] (97.52,116.94) .. controls (97.52,114.62) and (99.39,112.75) .. (101.7,112.75) .. controls (104.01,112.75) and (105.88,114.62) .. (105.88,116.94) .. controls (105.88,119.25) and (104.01,121.13) .. (101.7,121.13) .. controls (99.39,121.13) and (97.52,119.25) .. (97.52,116.94) -- cycle ;
	%Shape: Ellipse [id:dp8169827096481999] 
	\draw  [fill={rgb, 255:red, 0; green, 0; blue, 0 }  ,fill opacity=1 ] (135.91,77.79) .. controls (135.91,75.47) and (137.78,73.6) .. (140.09,73.6) .. controls (142.4,73.6) and (144.27,75.47) .. (144.27,77.79) .. controls (144.27,80.11) and (142.4,81.98) .. (140.09,81.98) .. controls (137.78,81.98) and (135.91,80.11) .. (135.91,77.79) -- cycle ;
	%Shape: Ellipse [id:dp9803689871178022] 
	\draw  [fill={rgb, 255:red, 0; green, 0; blue, 0 }  ,fill opacity=1 ] (136.89,116.94) .. controls (136.89,114.62) and (138.76,112.75) .. (141.07,112.75) .. controls (143.38,112.75) and (145.25,114.62) .. (145.25,116.94) .. controls (145.25,119.25) and (143.38,121.13) .. (141.07,121.13) .. controls (138.76,121.13) and (136.89,119.25) .. (136.89,116.94) -- cycle ;
	%Shape: Ellipse [id:dp09162610855419095] 
	\draw  [fill={rgb, 255:red, 0; green, 0; blue, 0 }  ,fill opacity=1 ] (175.28,116.61) .. controls (175.28,114.3) and (177.15,112.42) .. (179.46,112.42) .. controls (181.76,112.42) and (183.64,114.3) .. (183.64,116.61) .. controls (183.64,118.93) and (181.76,120.8) .. (179.46,120.8) .. controls (177.15,120.8) and (175.28,118.93) .. (175.28,116.61) -- cycle ;
	%Straight Lines [id:da776313181604554] 
	\draw    (102.35,38.97) -- (101.7,116.94) ;
	%Straight Lines [id:da21256404137119744] 
	\draw    (101.7,116.94) -- (179.46,116.61) ;
	%Shape: Ellipse [id:dp6817564370772887] 
	\draw  [fill={rgb, 255:red, 0; green, 0; blue, 0 }  ,fill opacity=1 ] (137.21,78.12) .. controls (137.21,75.8) and (139.08,73.92) .. (141.39,73.92) .. controls (143.7,73.92) and (145.57,75.8) .. (145.57,78.12) .. controls (145.57,80.43) and (143.7,82.31) .. (141.39,82.31) .. controls (139.08,82.31) and (137.21,80.43) .. (137.21,78.12) -- cycle ;
	%Shape: Ellipse [id:dp8424960424013599] 
	\draw  [fill={rgb, 255:red, 0; green, 0; blue, 0 }  ,fill opacity=1 ] (136.89,117.26) .. controls (136.89,114.95) and (138.76,113.07) .. (141.07,113.07) .. controls (143.38,113.07) and (145.25,114.95) .. (145.25,117.26) .. controls (145.25,119.58) and (143.38,121.46) .. (141.07,121.46) .. controls (138.76,121.46) and (136.89,119.58) .. (136.89,117.26) -- cycle ;
	%Shape: Ellipse [id:dp8314035611108691] 
	\draw  [fill={rgb, 255:red, 0; green, 0; blue, 0 }  ,fill opacity=1 ] (174.95,116.94) .. controls (174.95,114.62) and (176.82,112.75) .. (179.13,112.75) .. controls (181.44,112.75) and (183.31,114.62) .. (183.31,116.94) .. controls (183.31,119.25) and (181.44,121.13) .. (179.13,121.13) .. controls (176.82,121.13) and (174.95,119.25) .. (174.95,116.94) -- cycle ;
	%Straight Lines [id:da09907149260319081] 
	\draw    (358.72,40.27) -- (435.82,117.92) ;
	%Shape: Ellipse [id:dp3080373044159854] 
	\draw  [fill={rgb, 255:red, 0; green, 0; blue, 0 }  ,fill opacity=1 ] (354.53,40.27) .. controls (354.53,37.96) and (356.41,36.08) .. (358.72,36.08) .. controls (361.02,36.08) and (362.9,37.96) .. (362.9,40.27) .. controls (362.9,42.59) and (361.02,44.47) .. (358.72,44.47) .. controls (356.41,44.47) and (354.53,42.59) .. (354.53,40.27) -- cycle ;
	%Shape: Ellipse [id:dp4428650846734634] 
	\draw  [fill={rgb, 255:red, 0; green, 0; blue, 0 }  ,fill opacity=1 ] (354.21,79.42) .. controls (354.21,77.11) and (356.08,75.23) .. (358.39,75.23) .. controls (360.7,75.23) and (362.57,77.11) .. (362.57,79.42) .. controls (362.57,81.74) and (360.7,83.61) .. (358.39,83.61) .. controls (356.08,83.61) and (354.21,81.74) .. (354.21,79.42) -- cycle ;
	%Shape: Ellipse [id:dp30611485219817247] 
	\draw  [fill={rgb, 255:red, 0; green, 0; blue, 0 }  ,fill opacity=1 ] (353.88,118.24) .. controls (353.88,115.93) and (355.76,114.05) .. (358.06,114.05) .. controls (360.37,114.05) and (362.24,115.93) .. (362.24,118.24) .. controls (362.24,120.56) and (360.37,122.43) .. (358.06,122.43) .. controls (355.76,122.43) and (353.88,120.56) .. (353.88,118.24) -- cycle ;
	%Shape: Ellipse [id:dp4072937053376946] 
	\draw  [fill={rgb, 255:red, 0; green, 0; blue, 0 }  ,fill opacity=1 ] (392.27,79.1) .. controls (392.27,76.78) and (394.15,74.9) .. (396.45,74.9) .. controls (398.76,74.9) and (400.63,76.78) .. (400.63,79.1) .. controls (400.63,81.41) and (398.76,83.29) .. (396.45,83.29) .. controls (394.15,83.29) and (392.27,81.41) .. (392.27,79.1) -- cycle ;
	%Shape: Ellipse [id:dp6098573802983804] 
	\draw  [fill={rgb, 255:red, 0; green, 0; blue, 0 }  ,fill opacity=1 ] (393.25,118.24) .. controls (393.25,115.93) and (395.12,114.05) .. (397.43,114.05) .. controls (399.74,114.05) and (401.61,115.93) .. (401.61,118.24) .. controls (401.61,120.56) and (399.74,122.43) .. (397.43,122.43) .. controls (395.12,122.43) and (393.25,120.56) .. (393.25,118.24) -- cycle ;
	%Shape: Ellipse [id:dp43459512540276] 
	\draw  [fill={rgb, 255:red, 0; green, 0; blue, 0 }  ,fill opacity=1 ] (431.64,117.92) .. controls (431.64,115.6) and (433.51,113.72) .. (435.82,113.72) .. controls (438.13,113.72) and (440,115.6) .. (440,117.92) .. controls (440,120.23) and (438.13,122.11) .. (435.82,122.11) .. controls (433.51,122.11) and (431.64,120.23) .. (431.64,117.92) -- cycle ;
	%Straight Lines [id:da3402388277882751] 
	\draw    (358.72,40.27) -- (358.06,118.24) ;
	%Straight Lines [id:da041298461105825135] 
	\draw    (358.06,118.24) -- (435.82,117.92) ;
	%Shape: Ellipse [id:dp6148201640967023] 
	\draw  [fill={rgb, 255:red, 0; green, 0; blue, 0 }  ,fill opacity=1 ] (393.57,79.42) .. controls (393.57,77.11) and (395.45,75.23) .. (397.76,75.23) .. controls (400.06,75.23) and (401.94,77.11) .. (401.94,79.42) .. controls (401.94,81.74) and (400.06,83.61) .. (397.76,83.61) .. controls (395.45,83.61) and (393.57,81.74) .. (393.57,79.42) -- cycle ;
	%Shape: Ellipse [id:dp9386939336051117] 
	\draw  [fill={rgb, 255:red, 0; green, 0; blue, 0 }  ,fill opacity=1 ] (393.25,118.57) .. controls (393.25,116.25) and (395.12,114.38) .. (397.43,114.38) .. controls (399.74,114.38) and (401.61,116.25) .. (401.61,118.57) .. controls (401.61,120.88) and (399.74,122.76) .. (397.43,122.76) .. controls (395.12,122.76) and (393.25,120.88) .. (393.25,118.57) -- cycle ;
	%Shape: Ellipse [id:dp7106354251899841] 
	\draw  [fill={rgb, 255:red, 0; green, 0; blue, 0 }  ,fill opacity=1 ] (431.31,118.24) .. controls (431.31,115.93) and (433.19,114.05) .. (435.49,114.05) .. controls (437.8,114.05) and (439.67,115.93) .. (439.67,118.24) .. controls (439.67,120.56) and (437.8,122.43) .. (435.49,122.43) .. controls (433.19,122.43) and (431.31,120.56) .. (431.31,118.24) -- cycle ;
	%Straight Lines [id:da3904803443612095] 
	\draw    (358.39,79.26) -- (397.27,79.1) ;
	%Straight Lines [id:da09098394786160224] 
	\draw    (397.27,79.1) -- (397.43,118.57) ;
	%Straight Lines [id:da5317229073223257] 
	\draw    (358.39,79.26) -- (397.43,118.57) ;
	
	% Text Node
	\draw (84.71,66.72) node [anchor=north west][inner sep=0.75pt]   [align=left] {$\displaystyle 1$};
	% Text Node
	\draw (148.4,66.87) node [anchor=north west][inner sep=0.75pt]   [align=left] {$\displaystyle 1$};
	% Text Node
	\draw (93.89,124.64) node [anchor=north west][inner sep=0.75pt]   [align=left] {$\displaystyle 0$};
	% Text Node
	\draw (96.41,16.66) node [anchor=north west][inner sep=0.75pt]   [align=left] {$\displaystyle 0$};
	% Text Node
	\draw (107.67,176.47) node [anchor=north west][inner sep=0.75pt]   [align=left] {the heights};
	% Text Node
	\draw (138.89,125.6) node [anchor=north west][inner sep=0.75pt]   [align=left] {$\displaystyle 1$};
	% Text Node
	\draw (181.13,124.13) node [anchor=north west][inner sep=0.75pt]   [align=left] {$\displaystyle 0$};
	% Text Node
	\draw (326.5,167.47) node [anchor=north west][inner sep=0.75pt]   [align=left] {the induced subdivision};
	% Text Node
	\draw (365.95,63.41) node [anchor=north west][inner sep=0.75pt]  [font=\tiny] [align=left] {$\displaystyle Q_{1}$};
	% Text Node
	\draw (364,103.21) node [anchor=north west][inner sep=0.75pt]  [font=\tiny] [align=left] {$\displaystyle Q_{3}$};
	% Text Node
	\draw (377.58,88.2) node [anchor=north west][inner sep=0.75pt]  [font=\tiny] [align=left] {$\displaystyle Q_{2}$};
	% Text Node
	\draw (406.94,103.21) node [anchor=north west][inner sep=0.75pt]  [font=\tiny] [align=left] {$\displaystyle Q_{4}$};

\end{tikzpicture}
	\end{center}
 
    \caption{On the left, the Newton polygon of the polynomial $f$ from Example \ref{ex-dual}. Indicated are the negatives of the valuations of coefficients of $f$. On the right, the subdivision induced by these heights via projection from above.}
    \label{fig:subdivision}
\end{figure}
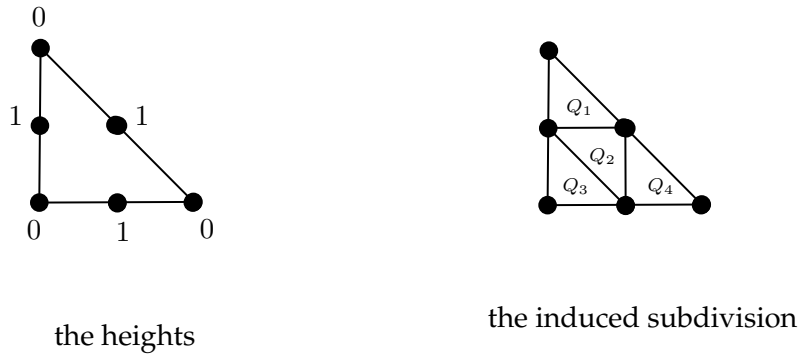

\end{example}

It is possible to construct two polynomials that lead to two different subdivisions of the same polygon, but for which the tropical curves coincide as subsets in $\RR^2$. This shows that it is inadequate to view tropical plane curves as subsets of $\RR^2$, but rather one should equip them with additional structure. First, we add integer weights on edges.

\begin{definition}
The \emph{weight} of an edge of a tropical plane curve equals the lattice length of the corresponding dual edge in the dual Newton subdivision. We call an edge \emph{even} if its weight is even.
\end{definition}

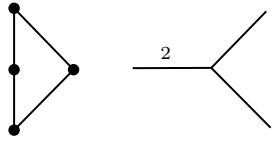
\begin{figure}
	\vspace{0.3cm}

	\begin{center}
		\tikzset{every picture/.style={line width=0.75pt}} %set default line width to 0.75pt        

\begin{tikzpicture}[x=0.75pt,y=0.75pt,yscale=-1,xscale=1]
	%uncomment if require: \path (0,784); %set diagram left start at 0, and has height of 784
	
	%Straight Lines [id:da6208660064350149] 
	\draw    (189.83,129.38) -- (217.53,100.85) ;
	%Straight Lines [id:da21256404137119744] 
	\draw    (150.2,129.52) -- (189.83,129.38) ;
	%Straight Lines [id:da09907149260319081] 
	\draw    (120.38,130.29) -- (90.38,99.29) ;
	%Straight Lines [id:da3402388277882751] 
	\draw    (189.83,129.38) -- (219.75,159.59) ;
	%Straight Lines [id:da041298461105825135] 
	\draw    (90.38,160.79) -- (90.38,99.29) ;
	%Straight Lines [id:da5317229073223257] 
	\draw    (120.38,130.29) -- (90.38,160.79) ;
	%Shape: Circle [id:dp5601260650531433] 
	\draw  [fill={rgb, 255:red, 0; green, 0; blue, 0 }  ,fill opacity=1 ] (88,99.29) .. controls (88,97.98) and (89.06,96.92) .. (90.38,96.92) .. controls (91.69,96.92) and (92.75,97.98) .. (92.75,99.29) .. controls (92.75,100.6) and (91.69,101.67) .. (90.38,101.67) .. controls (89.06,101.67) and (88,100.6) .. (88,99.29) -- cycle ;
	%Shape: Circle [id:dp3412741922932213] 
	\draw  [fill={rgb, 255:red, 0; green, 0; blue, 0 }  ,fill opacity=1 ] (88,130.29) .. controls (88,128.98) and (89.06,127.92) .. (90.38,127.92) .. controls (91.69,127.92) and (92.75,128.98) .. (92.75,130.29) .. controls (92.75,131.6) and (91.69,132.67) .. (90.38,132.67) .. controls (89.06,132.67) and (88,131.6) .. (88,130.29) -- cycle ;
	%Shape: Circle [id:dp35034359421934635] 
	\draw  [fill={rgb, 255:red, 0; green, 0; blue, 0 }  ,fill opacity=1 ] (88,160.79) .. controls (88,159.48) and (89.06,158.42) .. (90.38,158.42) .. controls (91.69,158.42) and (92.75,159.48) .. (92.75,160.79) .. controls (92.75,162.1) and (91.69,163.17) .. (90.38,163.17) .. controls (89.06,163.17) and (88,162.1) .. (88,160.79) -- cycle ;
	%Shape: Circle [id:dp5308501346824234] 
	\draw  [fill={rgb, 255:red, 0; green, 0; blue, 0 }  ,fill opacity=1 ] (118,130.29) .. controls (118,128.98) and (119.06,127.92) .. (120.38,127.92) .. controls (121.69,127.92) and (122.75,128.98) .. (122.75,130.29) .. controls (122.75,131.6) and (121.69,132.67) .. (120.38,132.67) .. controls (119.06,132.67) and (118,131.6) .. (118,130.29) -- cycle ;
	
	% Text Node
	\draw (162.7,117.27) node [anchor=north west][inner sep=0.75pt]  [font=\tiny] [align=left] {$\displaystyle 2$};

\end{tikzpicture}
	\end{center}

	\caption{On the right, a tropical plane curve having an edge of weight $2$. If we do not mention weights in a picture (as for the remaining two edges here), they are supposed to be one. On the left the dual subdivision having an edge of integer length $2$ dual to the edge of weight $2$.}
    \label{fig:weight}
\end{figure}

For an example, see Figure \ref{fig:weight}. With these weights, one discovers a beautiful feature of tropical plane curves, namely that they satisfy the following balancing condition. Let $v$ be a vertex in a tropical plane curve $C$. For each edge $e$ incident to $v$, let $w_e$ denote its weight and $p_e$ the primitive integral vector pointing from $v$ into the direction of $e$. Then
\begin{equation} \label{eq-balancing}
	\sum_{e \text{ incident to } v} w_e p_e=0.
\end{equation}

\begin{comment}

\begin{theorem}[Balancing condition]\label{thm-balancing}
	Let $C$ be a tropical plane curve and $v\in C$ a vertex. For each edge $e$ incident to $v$, let $w_e$ denote its weight and $p_e$ the primitive integral vector pointing from $v$ into the direction of $e$. Then
	$$\sum_{v \in \delta e} w_e p_e=0.$$
\end{theorem}

\end{comment}

The following example shows that we need to add even more structure to the datum of a tropical curve if we want to keep track of the genus of an algebraic curve: we should view tropical curves as being parametrized by an abstract graph, i.e. abstract tropical curves. 

\begin{example}\label{ex-tropsingcubic}
	Consider the cubic polynomial $p$ over the Puiseux series:
	$$p(x,y)=t^5x^3+xy^2+t^3y^3-xy-y^2+(-3t^5)\cdot x+(1-3t^3)\cdot y+(2t^3+2t^5) \in \Puiseux{\CC}[x,y].$$
	Figure \ref{fig:tropsingcubic} shows the tropicalization of the cubic.

	\begin{figure}
	   	\begin{center}
			\tikzset{every picture/.style={line width=0.75pt}} %set default line width to 0.75pt        

\begin{tikzpicture}[x=0.75pt,y=0.75pt,yscale=-1,xscale=1]
	%uncomment if require: \path (0,784); %set diagram left start at 0, and has height of 784
	
	%Straight Lines [id:da6208660064350149] 
	\draw    (189.83,129.38) -- (217.53,100.85) ;
	%Straight Lines [id:da21256404137119744] 
	\draw    (111.14,130.04) -- (189.83,129.38) ;
	%Straight Lines [id:da09907149260319081] 
	\draw    (110.86,70.9) -- (139.14,70.9) ;
	%Straight Lines [id:da3402388277882751] 
	\draw    (189.83,129.38) -- (147.71,207.75) ;
	%Straight Lines [id:da041298461105825135] 
	\draw    (139.8,200.22) -- (139.14,70.9) ;
	%Straight Lines [id:da5317229073223257] 
	\draw    (169.14,40.4) -- (139.14,70.9) ;
	%Straight Lines [id:da06714871372847098] 
	\draw    (111.52,200.22) -- (139.8,200.22) ;
	%Straight Lines [id:da6991592040747803] 
	\draw    (139.8,200.22) -- (147.71,207.75) ;
	%Straight Lines [id:da7062910910143697] 
	\draw    (148,234.04) -- (147.71,207.75) ;
	%Shape: Circle [id:dp048387054384323336] 
	\draw  [fill={rgb, 255:red, 0; green, 0; blue, 0 }  ,fill opacity=1 ] (-1.33,89.63) .. controls (-1.33,88.71) and (-0.59,87.97) .. (0.33,87.97) .. controls (1.25,87.97) and (2,88.71) .. (2,89.63) .. controls (2,90.55) and (1.25,91.3) .. (0.33,91.3) .. controls (-0.59,91.3) and (-1.33,90.55) .. (-1.33,89.63) -- cycle ;
	%Shape: Circle [id:dp27063879852987005] 
	\draw  [fill={rgb, 255:red, 0; green, 0; blue, 0 }  ,fill opacity=1 ] (-1.33,110.3) .. controls (-1.33,109.38) and (-0.59,108.63) .. (0.33,108.63) .. controls (1.25,108.63) and (2,109.38) .. (2,110.3) .. controls (2,111.22) and (1.25,111.97) .. (0.33,111.97) .. controls (-0.59,111.97) and (-1.33,111.22) .. (-1.33,110.3) -- cycle ;
	%Shape: Circle [id:dp8741842471116875] 
	\draw  [fill={rgb, 255:red, 0; green, 0; blue, 0 }  ,fill opacity=1 ] (-1,129.63) .. controls (-1,128.71) and (-0.25,127.97) .. (0.67,127.97) .. controls (1.59,127.97) and (2.33,128.71) .. (2.33,129.63) .. controls (2.33,130.55) and (1.59,131.3) .. (0.67,131.3) .. controls (-0.25,131.3) and (-1,130.55) .. (-1,129.63) -- cycle ;
	%Shape: Circle [id:dp5353930380131918] 
	\draw  [fill={rgb, 255:red, 0; green, 0; blue, 0 }  ,fill opacity=1 ] (-1.33,150.3) .. controls (-1.33,149.38) and (-0.59,148.63) .. (0.33,148.63) .. controls (1.25,148.63) and (2,149.38) .. (2,150.3) .. controls (2,151.22) and (1.25,151.97) .. (0.33,151.97) .. controls (-0.59,151.97) and (-1.33,151.22) .. (-1.33,150.3) -- cycle ;
	%Shape: Circle [id:dp3615718869331195] 
	\draw  [fill={rgb, 255:red, 0; green, 0; blue, 0 }  ,fill opacity=1 ] (18.33,109.97) .. controls (18.33,109.05) and (19.08,108.3) .. (20,108.3) .. controls (20.92,108.3) and (21.67,109.05) .. (21.67,109.97) .. controls (21.67,110.89) and (20.92,111.63) .. (20,111.63) .. controls (19.08,111.63) and (18.33,110.89) .. (18.33,109.97) -- cycle ;
	%Shape: Circle [id:dp05437636482865582] 
	\draw  [fill={rgb, 255:red, 0; green, 0; blue, 0 }  ,fill opacity=1 ] (18.44,129.86) .. controls (18.44,128.94) and (19.19,128.19) .. (20.11,128.19) .. controls (21.03,128.19) and (21.78,128.94) .. (21.78,129.86) .. controls (21.78,130.78) and (21.03,131.52) .. (20.11,131.52) .. controls (19.19,131.52) and (18.44,130.78) .. (18.44,129.86) -- cycle ;
	%Shape: Circle [id:dp7402315211641278] 
	\draw  [fill={rgb, 255:red, 255; green, 255; blue, 255 }  ,fill opacity=1 ] (18.67,149.97) .. controls (18.67,149.05) and (19.41,148.3) .. (20.33,148.3) .. controls (21.25,148.3) and (22,149.05) .. (22,149.97) .. controls (22,150.89) and (21.25,151.63) .. (20.33,151.63) .. controls (19.41,151.63) and (18.67,150.89) .. (18.67,149.97) -- cycle ;
	%Shape: Circle [id:dp3155140283504593] 
	\draw  [fill={rgb, 255:red, 255; green, 255; blue, 255 }  ,fill opacity=1 ] (38.89,130.08) .. controls (38.89,129.16) and (39.64,128.41) .. (40.56,128.41) .. controls (41.48,128.41) and (42.22,129.16) .. (42.22,130.08) .. controls (42.22,131) and (41.48,131.74) .. (40.56,131.74) .. controls (39.64,131.74) and (38.89,131) .. (38.89,130.08) -- cycle ;
	%Shape: Circle [id:dp15529782841538864] 
	\draw  [fill={rgb, 255:red, 255; green, 255; blue, 255 }  ,fill opacity=1 ] (38.67,149.97) .. controls (38.67,149.05) and (39.41,148.3) .. (40.33,148.3) .. controls (41.25,148.3) and (42,149.05) .. (42,149.97) .. controls (42,150.89) and (41.25,151.63) .. (40.33,151.63) .. controls (39.41,151.63) and (38.67,150.89) .. (38.67,149.97) -- cycle ;
	%Shape: Circle [id:dp035581447824303836] 
	\draw  [fill={rgb, 255:red, 0; green, 0; blue, 0 }  ,fill opacity=1 ] (58.67,149.97) .. controls (58.67,149.05) and (59.41,148.3) .. (60.33,148.3) .. controls (61.25,148.3) and (62,149.05) .. (62,149.97) .. controls (62,150.89) and (61.25,151.63) .. (60.33,151.63) .. controls (59.41,151.63) and (58.67,150.89) .. (58.67,149.97) -- cycle ;
	%Straight Lines [id:da2590680779378176] 
	\draw    (0.33,89.63) -- (60.33,149.97) ;
	%Straight Lines [id:da9271010350422306] 
	\draw    (0.33,89.63) -- (0.33,150.3) ;
	%Straight Lines [id:da8700857952753496] 
	\draw    (0.33,150.3) -- (40.34,150.08) -- (60.33,149.97) ;
	%Straight Lines [id:da4555363175676831] 
	\draw    (0.67,129.63) -- (20.11,129.86) ;
	%Straight Lines [id:da05033836494243049] 
	\draw    (0.33,110.3) -- (20,109.97) ;
	%Straight Lines [id:da6581946644610628] 
	\draw    (20,109.97) -- (20.11,129.86) ;
	%Straight Lines [id:da7580707471749853] 
	\draw    (20.11,129.86) -- (60.33,149.97) ;
	%Straight Lines [id:da5922516814362112] 
	\draw    (0.33,150.3) -- (20.11,129.86) ;
	
	% Text Node
	\draw (214.7,107.56) node [anchor=north west][inner sep=0.75pt]  [font=\tiny] [align=left] {$\displaystyle 2$};
	% Text Node
	\draw (137.67,227.37) node [anchor=north west][inner sep=0.75pt]  [font=\tiny] [align=left] {3};
	% Text Node
	\draw (102.43,115.18) node [anchor=north west][inner sep=0.75pt]  [font=\tiny] [align=left] {$\displaystyle ( 0,0)$};
	% Text Node
	\draw (103.19,182.51) node [anchor=north west][inner sep=0.75pt]  [font=\tiny] [align=left] {$\displaystyle ( 0,-3)$};
	% Text Node
	\draw (104.9,54.79) node [anchor=north west][inner sep=0.75pt]  [font=\tiny] [align=left] {$\displaystyle ( 0,3)$};
	% Text Node
	\draw (191.83,128.38) node [anchor=north west][inner sep=0.75pt]  [font=\tiny] [align=left] {$\displaystyle \left(\frac{5}{2} ,0\right)$};
	% Text Node
	\draw (148.79,198.38) node [anchor=north west][inner sep=0.75pt]  [font=\tiny] [align=left] {$\displaystyle \left(\frac{2}{3} ,-\frac{11}{3}\right)$};
	% Text Node
	\draw (-6.89,75.56) node [anchor=north west][inner sep=0.75pt]  [font=\tiny] [align=left] {$\displaystyle -3$};
	% Text Node
	\draw (-8,156.45) node [anchor=north west][inner sep=0.75pt]  [font=\tiny] [align=left] {$\displaystyle -3$};
	% Text Node
	\draw (-10.89,106.36) node [anchor=north west][inner sep=0.75pt]  [font=\tiny] [align=left] {$\displaystyle 0$};
	% Text Node
	\draw (-10.22,126.36) node [anchor=north west][inner sep=0.75pt]  [font=\tiny] [align=left] {$\displaystyle 0$};
	% Text Node
	\draw (22.67,101.24) node [anchor=north west][inner sep=0.75pt]  [font=\tiny] [align=left] {$\displaystyle 0$};
	% Text Node
	\draw (22.44,122.13) node [anchor=north west][inner sep=0.75pt]  [font=\tiny] [align=left] {$\displaystyle 0$};
	% Text Node
	\draw (52.89,155.02) node [anchor=north west][inner sep=0.75pt]  [font=\tiny] [align=left] {$\displaystyle -5$};
	% Text Node
	\draw (12.44,156.36) node [anchor=north west][inner sep=0.75pt]  [font=\tiny] [align=left] {$\displaystyle -5$};

\end{tikzpicture}
		\end{center}
	
	    \caption{On the right: the tropicalization of the cubic given by the polynomial $p$ in Example \ref{ex-tropsingcubic}. Indicated are edge weights and the coordinates of the vertices. On the left: its dual subdivision. }
	    \label{fig:tropsingcubic}
	\end{figure}
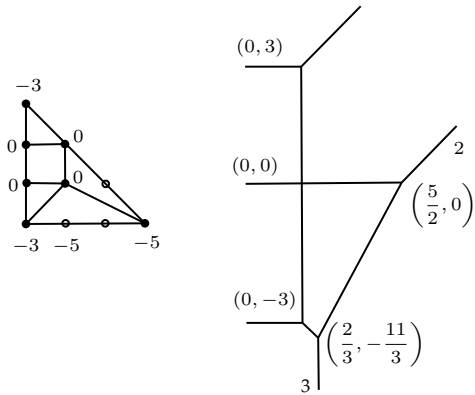
	
	The algebraic cubic curve $V(p)$ has a singularity at the point $(1,1)$
	%, as $$p(1,1)=\frac{\partial p}{\partial x}(1,1)=\frac{\partial p}{\partial y}(1,1)=0.$$ Thus, 
	and therefore $V(p)$ is rational.
	On the other hand, the tropicalization seems to have a cycle, suggesting genus 1. But we should  interpret the vertex dual to the square as a crossing of two edges, not as a vertex. We should parametrize our plane tropical curves by abstract graphs. If we do this, we can use a rational abstract graph, a tree, to parametrize the tropical cubic above.
\end{example}

%The example above motivates why we should think of tropical plane curves as being parametrized by an abstract graph, in order to be able to keep track of genus under tropicalization.
%The abstract graphs are viewed as abtract tropical curves.

When the lengths of the edges of a tropical curve are varied in a family, cycles of graphs can vanish into vertices. The correct definition of an abstract tropical curve therefore includes a genus function on the vertices keeping track of such issues. For our purposes, it is sufficient to restrict to abstract tropical curves with all vertex genus equal to 0. Such abstract tropical curves are called \emph{explicit}.

\begin{definition}%[Abstract tropical curves]
	An (explicit) \emph{abstract tropical curve} is a metric graph $\Gamma$ with unbounded edges of infinite length, called \emph{ends}. The \emph{genus} of a (connected) tropical curve is the first Betti number of the underlying graph, i.e. for connected $\Gamma$
	$$ g(\Gamma) = 1-\#\mbox{vertices}+\#\mbox{bounded edges}$$
	and if $\Gamma$ is disconnected with connected components $\Gamma_1,\ldots,\Gamma_r$ of genera $g_i := g(\Gamma_i)$, then the genus of $\Gamma$ is defined to be 
	$$g(\Gamma):= \sum_{i=1}^rg_i-r+1.$$
\end{definition}

\begin{definition}%[Parametrized tropical plane curve/tropical stable map]
A \emph{parametrized tropical plane curve} (also called a tropical stable map to $\RR^2$) is a tuple $(\Gamma,\varphi)$, such that $\Gamma$ is an abstract tropical curve and $\varphi$ is a map which is locally on each edge integer affine linear such that the balancing condition~\eqref{eq-balancing} is satisfied at every vertex. The weights are given here as expansion factors, i.e.\  an edge $e$ of length $l$ is mapped to a segment connecting a point $a\in \RR^2$ with $a+l\cdot w_e\cdot p_e$, where $p_e$ is the primitive integer vector pointing in the direction of $\varphi(e)$.

The \emph{genus} of $(\Gamma,\varphi)$ is defined to be the genus of $\Gamma$.

When we speak of \emph{vertices} of a parametrized tropical curve, we mean the vertices of the underlying abstract tropical curve.
\end{definition}

If we parametrize the cubic in Example~\ref{ex-tropsingcubic} with a graph, it is not clear whether the down end should be parametrized by three ends of weight 1 each or one end of weight 3. To specify this, we work with the slightly refined notion of Newton fans rather than their (roughly) dual Newton polygons when working with parametrized tropical curves.

\begin{definition}
A \emph{Newton fan} is a multiset $\Delta=\{v_1,\ldots,v_k\}$ of vectors $v_i\in \ZZ^2$ satisfying $\sum_{i=1}^k v_i=0$. 
We use the notation $\Delta=\{v_1^{m_1},\ldots, v_k^{m_k}\}$ to indicate
that the vector $v_i$ appears $m_i$ times in $\Delta$.
\end{definition}

If $\Delta= \{(a_1,b_1),\ldots,(a_k,b_k)\}$ is a Newton fan of vectors in $\ZZ^2$, one can
construct the \emph{dual polygon} $\Pi_\Delta$ in $\RR^2$ in the following
way: for each primitive integer direction $(a,b)$ in $\Delta$,
consider the vector $w\cdot (-b,a)$, where $w$ is the sum of the
weights of all vectors in $\Delta$ with primitive integer direction
$(a,b)$.
Now $\Pi_\Delta$ is the unique (up to translation) polygon whose
oriented edges (the orientation is induced by the usual orientation of
$\RR^2$) are exactly the vectors $w\cdot (-b,a)$. 

\begin{definition}
	The \emph{degree} of a parametrized tropical curve $(\Gamma,\varphi)$ is the Newton fan of directions of its ends. 
	
	 We say that $(\Gamma, \varphi)$ has \emph{degree $d$} if this multiset consists of $d$ times the vector $(-1,0)$, $d$ times $(0,-1)$ and $d$ times $(1,1)$. 
	These are tropical curves that arise as tropicalizations of degree $d$ curves in the projective plane $\mathbb{P}^2$ which intersect the toric boundary generically.
	
	The \emph {combinatorial type} of $(\Gamma, \varphi)$ is obtained by forgetting the metric information, i.e. we keep the topological type of $\Gamma$ together with the degree. 
\end{definition}

In our statements of correspondence theorems, by abuse of notation, we use $\Delta$ to denote a Newton fan for counts of tropical curves as well as for the curve class in the toric surface defined by the polygon $\Pi_\Delta$ via hyperplane sections. As correspondence theorems only cover cases where all ends have weight $1$, no confusion should occur.

\subsection{Complex and real counts of tropical plane curves}\label{subsec-counts}

Fix a Newton fan $\Delta$ and some genus $g$.
Let $n=\# \Delta +g-1$ and fix $n$ points $p_1,\ldots,p_n$ in tropical general position in $\mathbb{R}^2$ (see Definition 4.7 in \cite{Mi03} resp.\ Definition 5.33 in \cite{Ma06}).

For some piecewise linear graphs in $\RR^2$, there are several ways to parametrize them by an abstract tropical curve. Simple tropical curves can uniquely be parametrized, once we fix the convention that every point dual to a parallelogram should be viewed as a crossing of two edges and not as a vertex.

\begin{definition}[Simple tropical plane curves]
A parametrized tropical plane curve $(\Gamma,\varphi)$ is called \emph{simple} if $\Gamma$ is $3$-valent and the Newton subdivision dual to $\varphi(\Gamma)\subset \RR^2$ contains only triangles and paralellograms. 
\end{definition}

Here, for a tropical plane curve to be simple we do not insist that all ends have weight 1, as it is the case in some places in the literature.

As the following proposition shows, restricting to simple tropical plane curves is sufficient for the purpose of counting curves satisfying point conditions.

\begin{proposition}\label{prop-genericsimple}
Fix $\#\Delta+g-1$ points in $\mathbb{R}^2$ in tropical general position. Then any tropical plane curve of degree $\Delta$ and genus $g$ passing through these points is simple.
\end{proposition}
For a proof, see e.g.\ Lemma 5.34 in \cite{Ma06}.

\begin{definition}

Consider a simple tropical plane curve $C$. We define the \emph{complex multiplicity of $C$}
to be
$$\mult_{\mathbb{C}}(C) = \prod_\delta \Area(\delta)$$
where the product runs over all triangles $\delta$ in the dual Newton subdivision of $C$ and $\Area(\delta)$ denotes the normalized (lattice) area of $\delta$.

We define the \emph{real multiplicity of $C$}
to be
$$\mult_{\mathbb{R}}(C) =\begin{cases}0 &\mbox{if $C$ has an even edge}\\ \prod_\delta (-1)^{i(\delta)}&\mbox{if all edges of $C$ are odd},\end{cases}$$
where the product runs over all triangles $\delta$ in the dual Newton subdivision of $C$ and $i(\delta)$ denotes the number of interior lattice points of the triangle $\delta$.
\end{definition}

Note that by Pick's formula, the existence of an even edge is equivalent to the complex multiplicity being even.

By Mikhalkin's correspondence theorem \cite{Mi03} we can recover algebro-geometric counts from the tropical point of view when using $\mult_\CC$ and $\mult_\RR$. More precisely, denote by $N(\Delta,g)$ the Gromov-Witten invariant given by the number of complex genus $g$ curves of class $\Delta$ in the toric surface defined by $\Delta$ passing through $n$ points in general position. Moreover, for a real plane curve $\mathcal{C}$ define $s(\mathcal{C})$ to be the number of solitary nodes, i.e.\ nodes which are locally of the form $V(x^2+y^2=0)$.
If $g = 0$ and $\Delta$ defines a toric del Pezzo surface, then we denote by $W(\Delta,0)$ the Welschinger invariant, counting rational real plane curves of class $\Delta$ satisfying $n$ real point conditions, where each curve is counted with sign $(-1)^{s(\mathcal{C})}$. Mikhalkin's theorem says:

\begin{theorem}[Correspondence Theorem, \cite{Mi03}]\label{thm-corres-old}
Fix a Newton fan $\Delta$ whose entries have all weight one and such that the corresponding polygon defines a smooth toric surface. Fix a genus $g$ and $n:=\#\Delta+g-1$ points in general position in $\mathbb{R}^2$.
Then
$$N^\trop(\Delta,g):=\sum_C \mult_\CC (C)= N(\Delta,g),$$
where the sum runs over all (necessarily simple) parametrized tropical plane curves of degree $\Delta$ and genus $g$ passing through the points.

Moreover, in genus $g = 0$ and for $\Delta$ corresponding to a smooth del Pezzo surface we have
$$W^\trop(\Delta,0):=\sum_C \mult_\RR (C)= W(\Delta,0),$$ 
where, as before, the sum runs over all (necessarily simple) parametrized tropical plane curves of degree $\Delta$ and genus $g$ passing through the points.

\end{theorem}

Note that a similar correspondence theorem for the real case  can also be formulated in higher genus, however, $W(\Delta,g)$ is not an invariant but depends on the chosen point conditions. The fact that the corresponding tropical number is invariant is not a contradiction: this merely reflects the fact that the signed count of real curves passing through point conditions close to the tropical limit (i.e.\ which tropicalize to points in tropical general position) is invariant.

\begin{example}\label{ex-ninecubics}
Figure \ref{fig-ninecubics} shows $8$ points in $\mathbb{R}^2$ in tropical general position, and the nine rational tropical cubics which pass through them. The one in the last row on the right has complex multiplicity $4$ and real multiplicity $0$, all others have complex and real multiplicity equal to $1$. We obtain $N^\trop(3,0)= 8\cdot 1 + 1 \cdot 4 =12$ which equals, as expected, the number of complex rational cubics passing through $8$ points. We also obtain $W^\trop(3,0)=8$ which equals the Welschinger invariant of rational cubics passing through $8$ real points.

\begin{figure}
\begin{center}
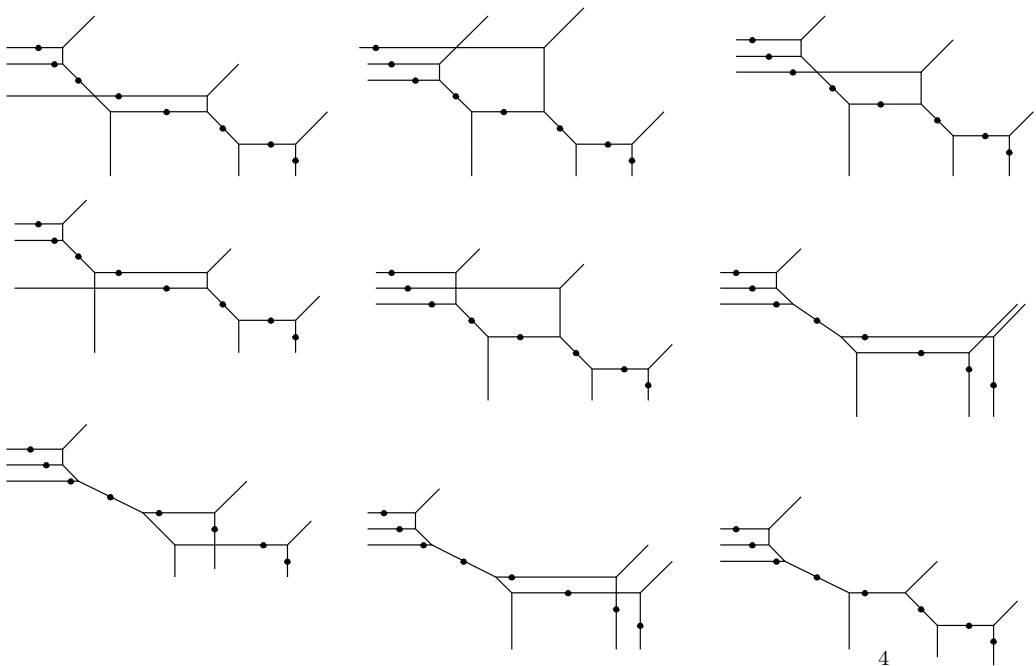
\end{center}
\caption{The picture shows all 9 simple rational tropical cubics through $8$ given points in general position. The one in the last row on the right has complex multiplicity $4$ and real multiplicity $0$. All others have complex and real multiplicity $1$.}\label{fig-ninecubics}
\end{figure}

\end{example}

\subsection{Arithmetic counts of tropical plane curves}

Let $K$ be a field.
We now recall the definition of the Grothendieck-Witt ring $\GW(K)$ in which arithmetic counts take their values.  

A \emph{quadratic space} is a finite-dimensional $K$-vector space $V$ equipped with a nondegenerate symmetric bilinear form $q \colon V \times V \to K$. Two quadratic spaces $(V,q)$ and $(V',q')$ are isomorphic if there is an isomorphism of $K$-vector spaces $\phi\colon V \to V'$ such that $q(v,w) = q'(\phi(v), \phi(w))$ for all $v$, $w$ in $V$. 
The set of isomorphism classes of quadratic spaces has the structure of a monoid with respect to direct sum $\oplus$.

\begin{definition}
The  \emph{Grothendieck-Witt ring} $\GW(K)$ is the Grothendieck group of the monoid of quadratic spaces described above. It is a commutative ring with the tensor product as multiplication.
\end{definition}

\begin{definition}
For $a \in K^\times$, we write $\langle a \rangle$ for the $1$-dimensional quadratic space $(K,q)$ with $q(x,y) = axy$.   The {\em hyperbolic plane} is $\mathbb{H} = \langle 1 \rangle + \langle -1 \rangle.$
\end{definition}

\begin{lemma}
\label{lem-GWrelations}
As an additive group, $\GW(K)$ is generated by $\big\{\langle a \rangle : a\in K^\times \big\}$, with relations generated by
\begin{enumerate}
\item $\langle a \rangle= \langle a b^2 \rangle$ for all $a$, $b$ in $K^\times$, and
\item $\langle a\rangle + \langle b \rangle= \langle a+b \rangle+ \langle ab(a+b) \rangle$,
for all $a$, $b$ in $K^\times$ such that $a + b \neq 0$.
\end{enumerate}
Furthermore, the multiplication satisfies the following relation.
\begin{enumerate}
\setcounter{enumi}{2}
\item $\langle a\rangle\langle b\rangle=\langle ab\rangle$ for $a,b\in K^\times$.
\end{enumerate}
\end{lemma}
For a proof, see \cite[Theorem 2.1.11, Remark 2.1.12]{deglise2023notes}.
 An easy consequence is:
\begin{lemma}\label{lem:hyperbolic} 
For any $a \in K^\times$, the quadratic space $\langle a \rangle +\langle -a \rangle$ is isomorphic to $\mathbb{H}$, as is the quadratic space $K^2$ with quadratic form $\big((x_1, x_2),(y_1, y_2)\big) \mapsto ax_1y_2 + ax_2y_1$.
\end{lemma}

%\begin{proof}
%The isomorphisms are induced by the change of variables with matrices presented in the equalities of products
 %   \[\frac{1}{2a}\begin{pmatrix}1+a&1-a\\1-a&1+a\end{pmatrix}\cdot \begin{pmatrix}a&0\\0&-a\end{pmatrix}\cdot \frac{1}{2a}\begin{pmatrix}1+a&1-a\\1-a&1+a\end{pmatrix}=
%    \begin{pmatrix}1&0\\0&-1\end{pmatrix}\] and \[ \frac{1}{2a}
%    \begin{pmatrix}2&a\\-2&a\end{pmatrix}\cdot \begin{pmatrix}0&a\\a&0\end{pmatrix}\cdot \frac{1}{2a} \begin{pmatrix}2&-2\\a&a\end{pmatrix}=
%    \begin{pmatrix}1&0\\0&-1\end{pmatrix}.\]   
%\end{proof}

\noindent When no confusion seems possible, we write $\langle a \rangle$, $\mathbb{H}$, and so on, not only for a given quadratic space, but also for its class in $\GW(K)$. 

Assume $K' / K$ is a finite separable field extension. Then any finite dimensional $K'$-vector space $V$ is also finite-dimensional as a $K$-vector space, and we write $V_K$ to denote $V$, viewed as a $K$-vector space. If $(V,q)$ is a quadratic space over $K'$, then $(V_K, \Tr_{K'|K} \circ q)$ is a quadratic space over $K$. One writes 
\[
\Tr_{K'|K} \colon \GW(K') \longrightarrow \GW(K)
\]
for the induced map of Grothendieck-Witt rings.  

\begin{example} \label{ex-GW-CC}
The Grothendieck-Witt ring of the complex numbers equals $\GW(\mathbb{C})\simeq \ZZ$ with generator $\langle 1 \rangle$ and the identification is given by evaluating the rank of a quadratic form. The Grothendieck-Witt ring of the reals equals $\GW(\mathbb{R})\simeq \ZZ[\mu_2]$, where $\mu_2$ denotes the cyclic group with $2$ elements. It is generated by $\langle 1\rangle$ and $\langle -1\rangle$.
\end{example}

%Methods from $\mathbb{A}^1$-homotopy allow to get meaningful results in enumerative geometry over an arbitrary base field $K$ when \enquote{counting} in $\GW(K)$ instead of the integers. To do this one assigns a \emph{quadratic weight} to each object one is counting. In the case of counting rational curves embedded in a toric del Pezzo surface, one gets an invariant count when counting each curve with the following quadratic enrichment of the Welschinger sign \cite{Wel05} defined by Levine \cite{LevineWelschinger}: 

Methods from $\mathbb{A}^1$-homotopy theory allow to define meaningful invariants for enumerative problems posed over an arbitrary base field $K$. In contrast to the real and complex situation, these arithmetic invariants live in $\GW(K)$ rather than the integers. For example, let $\mathcal{C}$ be a rational curve defined over $K$ (where $K$ is a perfect field with $\chara(K)\neq 2,3$) and embedded in a toric del Pezzo surface.  In this case the following quadratic enrichment of the Welschinger sign defined by Levine \cite{LevineWelschinger} is such an invariant:

\[\operatorname{Wel}_{\mathbb{A}^1}(\mathcal{C})= \Big\langle\prod_{\text{nodes } z}N_{\kappa(z)/K} \big(-\det \operatorname{Hessian}f(z) \big) \Big\rangle \in \GW(K)\]
where $\kappa(z)$ is the residue field of the node $z$ and $N_{\kappa(z)/K}$ the field norm. This naturally generalizes the Welschinger sign $\operatorname{Wel}(\mathcal{C})$: if $\mathcal{C}$ is defined over $\RR$, then $\langle \operatorname{Wel}(\mathcal{C})\rangle=\operatorname{Wel}_{\mathbb{A}^1}(\mathcal{C})$. Now define the quadratically weighted count of rational curves on a smooth toric del Pezzo surface as 
\begin{equation} \label{eq:arithmetic_count}
	N_{\mathbb{A}^1}(\Delta,0):=\sum_{\mathcal{C}} \operatorname{tr}_{\kappa(\mathcal{C})|K}\big(\operatorname{Wel}_{\mathbb{A}^1}(\mathcal{C})\big) \in \GW(K),
\end{equation}
where $\kappa(\mathcal{C})$ is the field of definition of $\mathcal{C}$ and the sum runs over all $\mathcal{C}$ which pass through the correct number of points in general position. It is shown in \cite{KassLevineSolomonWickelgren, KLSWOrientation} that $N_{\AA^1}(\Delta, 0)$ is independent of the choice of points.

The fact that $\GW(\Puiseux{K}) \simeq \GW(K)$ for $\chara(K)\neq 2$ (see Theorem 4.7 in \cite{MPS22}) hints to the fact that tropicalization and quadratic enrichment should \enquote{work well together}. Case studies of this fruitful relation have already appeared in \cite{puentes2022quadratically} and \cite{MPS22}. Recently, the first and third author have also proved that this relation can be used in the enumerative context, by proving a quadractically enriched correspondence theorem which specializes to Theorem \ref{thm-corres-old} for the field of complex resp.\ real numbers. 
The tropical curves have to be counted with an arithmetic multiplicity in this context, which we now define.

\begin{definition}%[Arithmetic multiplicity of a simple tropical plane curve]
	\label{def-troparith}
Consider a simple tropical plane curve $C$. Assume the ends of $C$ have all weight $1$. 
We define the \emph{arithmetic multiplicity of $C$}
to be
$$\mult_{\mathbb{A}^1}(C) = \begin{cases} \frac{\mult_{\mathbb{C}}(C)-1}{2}\cdot \mathbb{H}+ \big\langle \mult_{\mathbb{R}}(C) \big\rangle & \mbox{ if $\mult_\CC(C)$ is odd,}\\
\frac{\mult_{\mathbb{C}}(C)}{2}\cdot \mathbb{H} & \mbox{ if $\mult_\CC(C)$ is even.}

 \end{cases}$$
\end{definition}

\begin{theorem}[Quadratically enriched correspondence theorem, \cite{JPP23}]
\label{thm-corres-new}

Fix a Newton fan $\Delta$ whose entries have all weight 1 and such that the corresponding polygon defines a smooth toric del Pezzo surface. 
 Let $K$ be a perfect field of large enough characteristic, i.e.\ larger than the diameter of $\Delta$.
Fix $n=\#\Delta-1$ points in general position in $\mathbb{R}^2$.
Then
$$N^\trop_{\mathbb{A}^1}(\Delta,0):=\sum_C \mult_{\mathbb{A}^1} (C)= N_{\mathbb{A}^1}(\Delta,0)\in \GW(K),$$
where the sum runs over all (necessarily simple) parametrized rational tropical plane curves of degree $\Delta$  passing through the points, and
where $N_{\mathbb{A}^1}(\Delta,0)$ denotes the arithmetic count defined in~\eqref{eq:arithmetic_count} of rational curves of class $\Delta$ in the toric surface defined by $\Delta$ passing through $n$ points in general position.
\end{theorem}
 The restriction on the charactertistics of the field is coming from the fact that we have contributions involving weights of edges in the arithmetic multiplicity, which should not be $0$. With the assumption that the characteristic is larger than the diameter of $\Delta$, no appearing edge weight can be $0$.

\begin{example}
Using the tropical plane curve count from Example \ref{ex-ninecubics}, Figure \ref{fig-ninecubics}, we obtain $$N^\trop_{\mathbb{A}^1}(3,0)=2\mathbb{H} + 8\cdot \langle 1 \rangle.$$
\end{example}

\section{Invariance of the arithmetic count of tropical plane curves}
Firstly, we generalize Definition \ref{def-troparith} to arbitrary Newton fans.

\begin{definition}%[Arithmetic multiplicity of a simple tropical plane curve]
	\label{def-troparith2}
Consider a simple tropical plane curve $C$  and $K$ a field. Assume the ends of $C$ have weights $w_1,\ldots,w_k$  and that the characteristic of $K$ does not divide any of the $w_i$.
We define the \emph{arithmetic multiplicity of $C$}
to be
$$\mult_{\mathbb{A}^1}(C) = 
\begin{cases} 
	\frac{\mult_{\mathbb{C}}(C)-1}{2}\cdot \mathbb{H}+ \big\langle (-1)^i  w_1 \cdots w_k \big\rangle & \mbox{ if $\mult_\CC(C)$ is odd,}\\
	\frac{\mult_{\mathbb{C}}(C)}{2}\cdot \mathbb{H} & \mbox{ if $\mult_\CC(C)$ is even,}
\end{cases}$$
where $i$ denotes the total number of interior lattice points appearing in triangles of the Newton subdivision dual to $C$.
\end{definition}

\begin{remark}
In the quadratically enriched correspondence theorem \cite{JPP23} the first and third author only consider tropical curves with ends of weight one. In this case the multiplicity in Definition~\ref{def-troparith2} agrees with Definition~\ref{def-troparith}. However, if one allowed tropical curves with higher odd weight ends in the proof of the correspondence theorem in \cite{JPP23}, then one would get the above multiplicity. 
In case of even weights, this is not necessarily the case. However, in this case Definition 3.1 is the natural generalization of the real multiplicity since in this case the signature equals $0$.
\end{remark}

The arithmetic multiplicity of a tropical curve can also be written as a product of vertex multiplicities.

\begin{lemma}\label{lem-vertexmult}
Let $C$ be a simple tropical plane curve and $V$ one of its vertices. By cutting the three edges adjacent to $V$ and thinking of them as unbounded ends, we can view a neighborhood of $V$ as a tropical curve itself, which, abusing notation, we also denote by $V$. Its arithmetic multiplicity $\mult_{\mathbb{A}^1}(V)$ is defined in Definition~\ref{def-troparith2}.
We have
$$ \mult_{\mathbb{A}^1}(C) = \prod_{V \text{ vertex of } C} \mult_{\mathbb{A}^1}(V).$$
\end{lemma}

\begin{proof}
We prove the lemma by induction on the number of vertices. For one vertex, the statement is clear. Now assume that $C$ has $k$ vertices. If $C$ has a bridge edge, we can cut it, producing two connected tropical curves $C_1$ and $C_2$ for each of which the induction assumption holds. The vertices of $C$ can be divided into vertices that belong to $C_1$ or $C_2$ after cutting. 
Assume first that the complex multiplicity of both $C_1$ and $C_2$ (and with this, the complex multiplicity of $C$) is odd. Assume that the ends of $C_1$ have weights $w_1\ldots,w_l,w$, where $w$ is the weight of the cut edge, and the ends of $C_2$ have weights $w, w_{l+1},\ldots,w_k$. Accordingly, the ends of $C$ have weights $w_1, \ldots,w_k$. Then we have 
\begin{multline} \label{eq:prod_of_multiplicities}
	\prod_{V\in C} \mult_{\mathbb{A}^1}(V) 
	= \prod_{V\in C_1} \mult_{\mathbb{A}^1}(V)\cdot \prod_{V\in C_2} \mult_{\mathbb{A}^1}(V) \\ 
	= \Big(\frac{\mult_{\mathbb{C}}(C_1)-1}{2}\cdot \mathbb{H}+ \Big\langle (-1)^{i_1} \cdot \prod_{j=1}^l w_j\cdot w \Big\rangle \Big) 
	\cdot \Big(\frac{\mult_{\mathbb{C}}(C_2)-1}{2}\cdot \mathbb{H}+ \Big\langle (-1)^{i_2} \cdot \prod_{j=l
+1}^k w_j\cdot w \Big\rangle \Big) ,
\end{multline}
where the second equality holds by the induction assumption. Here, for $m={1,2}$, we use $i_m$ to denote the number of interior points in triangles dual to vertices of $C_m$.
Multiplication with $\mathbb{H}$ in the Grothendieck-Witt ring gives an integer multiple of $\mathbb{H}$ (apply (3) in Lemma \ref{lem-GWrelations} and Lemma \ref{lem:hyperbolic}), thus expanding the right hand side of Equation~\eqref{eq:prod_of_multiplicities} will give an expression of the form 
\begin{equation} \label{eq-exprGW}
	M\cdot \mathbb{H}+ \Big\langle (-1)^{i_1+i_2} \cdot \prod_{j=1}^k w_j \Big\rangle ,
\end{equation}
where $M$ is an integer. Note that in the second summand, a factor of $w^2$ can be canceled because of Lemma~\ref{lem-GWrelations}. Furthermore, the sum of the numbers of interior points in triangles dual to vertices of $C_1$ and of $C_2$ obviously equals the number of interior points in triangles dual to vertices of $C$.

It remains to compute $M$. 
We can do so, by specializing to $K = \CC$ (which is equivalent to computing ranks of the quadratic forms, see Example~\ref{ex-GW-CC}). For the right hand side of expression~\eqref{eq:prod_of_multiplicities} this yields
$$ \Big(\frac{\mult_{\mathbb{C}}(C_1)-1}{2}\cdot 2+1 \Big) \cdot \Big(\frac{\mult_{\mathbb{C}}(C_2)-1}{2}\cdot 2+1\Big) = \mult_{\mathbb{C}}(C_1)\cdot \mult_{\mathbb{C}}(C_2) = \mult_{\mathbb{C}}(C) $$
and for expression~\eqref{eq-exprGW} we obtain
$ M\cdot 2+1$. As these two must be equal, we deduce $M=\frac{ \mult_{\mathbb{C}}(C)-1}{2}$ and the claim follows.

The cases where one or both of the $C_i$ have an even complex multiplicity follows analogously.
If $C$ does not have a bridge edge, we have to cut at least two edges to produce two connected components. The argument follows along the same lines, only in expression \eqref{eq-exprGW} we cancel not only one contribution of the form $w^2$, where $w$ is the weight of a cut edge, but two or more such contributions. 
\end{proof}

\begin{definition}[Arithmetic count of plane tropical curves]
Fix a Newton fan $\Delta$ and a genus $g$.  Let $K$ be a field of characteristic larger than the diameter of $\Delta$. Fix $n=\#\Delta+g-1$ points in general position in $\mathbb{R}^2$.
Then
$$N^\trop_{\mathbb{A}^1}(\Delta,g):=\sum_C \mult_{\mathbb{A}^1} (C)\in \GW(K),$$
where the sum runs over all (necessarily simple) parametrized tropical plane curves $C$ of degree $\Delta$ and genus $g$ passing through the points.
\end{definition}

Just as in the real case, these numbers (for the case of a Newton fan with ends of weight $1$ only) could also appear in a correspondence theorem: the arithmetic count is not an invariant, but if we arithmetically count curves through points close to the tropical limit, we obtain the tropical count.

Theorem \ref{thm-invariance} shows that the arithmetic count of tropical curves of degree $\Delta$ and genus $g$ is well-defined, i.e.\ it does not depend on the position of the point conditions.

\begin{remark}[Arithmetic count of tropical plane curves via lattice path]
As a consequence of Lemma \ref{lem-vertexmult} and Theorem~\ref{thm-invariance}, we can use Mikhalkin's lattice path algorithm (see Section 7 of \cite{Mi03}) to compute the number $N^\trop_{\mathbb{A}^1}(\Delta,g)$ for any $\Delta$ and $g$: we just have to use the arithmetic multiplicity of a triangle whenever the area (i.e. the complex multiplicity) occurs  in Mikhalkin's formulation. For an example see Corollary~\ref{cor:lattice_path}.
% instead of its area (as complex multiplicity) resp.\ $(-1)^i$ where $i$ is the number of its interior points (as real multiplicity).
\end{remark}

\begin{corollary} \label{cor:lattice_path}
For the arithmetic count of tropical curves of degree $d$ of maximal genus (i.e.\ smooth curves) resp.\ of genus one less than maximal (i.e.\ curves with one node) we obtain
\begin{align*}
	N^\trop_{\mathbb{A}^1}\Big(d,\binom{d-1}{2}\Big) &=\langle 1\rangle, \quad \mbox{ and }\\  N^\trop_{\mathbb{A}^1}\Big(d,\binom{d-1}{2}-1\Big) &= (d-1) (d-2)\cdot \mathbb{H}+(d^2-1)\cdot \langle 1 \rangle. \end{align*}
\end{corollary}

\begin{proof}
This follows from a lattice path count. We are working in the triangle $\Delta_d$ with vertices $(0,0)$, $(d,0)$ and $(0,d)$. In the first case, we have so many point conditions that the corresponding lattice path has to use every point in $\Delta_d$. There is a unique such path, and every triangle that appears in the unique underlying Newton subdivision is of area $1$. 

In the second case, the lattice path must miss exactly one point. If it misses an interior point, we get two triangles of weight $2$ each, leading to a contribution of $2\mathbb{H}$ for each interior point of which there are $\binom{d-1}{2}$. Otherwise, we can miss a point on the bottom edge or on the diagonal edge. In each case, we obtain one parallelogram in the subdivision and else only triangles of area $1$. Each such subdivision contributes $\langle 1 \rangle$. If we miss the point $(i,0)$ (where $i$ ranges from $0$ to $d-2$), then there are $d-i-1$ ways to fit in a parallelogram. If we miss the point $(d-i,i)$ (where $i$ ranges from $1$ to $d-1$), then there are $d-i+1$ ways to fit in a parallelogram. Adding up the contributions, we obtain the result.
\end{proof}

We now state the main theorem of this section:

\begin{theorem}[Invariance of arithmetic counts of plane tropical curves]\label{thm-invariance}
 Let $K$ be a field of large enough characteristic, i.e.\ larger than the diameter of $\Delta$.
The count $N^\trop_{\mathbb{A}^1}(\Delta,g)$ does not depend on the location of the points in $\RR^2$ (as long as they are in general position) through which we require the tropical curves to pass.
\end{theorem}

We emphasize that the combinatorial side of the tropical plane curve count is the same for a complex, a real, and an arithmetic count. This allows us to use already established techniques for tropical plane curve counts in the proof of Theorem~\ref{thm-invariance}. More concretely, in Section~3 \cite{GM051}, a moduli space $\overline{\mathcal{M}}_{g, \Delta}$ of plane tropical curves with marked points is constructed and with it evaluation maps
\[ \operatorname{ev} : \overline{\mathcal{M}}_{g, \Delta} \longrightarrow (\RR^2)^n \]
which evaluate the position of the marked points in $\mathbb{R}^2$. 
Roughly, this moduli spaces is a union of polyhedra, where the union runs over combinatorial types (here, in addition to the underlying graph and all directions for the edges, we also remember the edges resp.\ vertices on which the marked points lie). Since the (complex, real or arithmetic) multiplicity only depends on the combinatorial type of a tropical plane curve, the function which counts inverse images of the evaluation map with their multiplicity  is locally constant.
In fact, if we fix points in general position, the tropical curves meeting them all live in the interior of a top-dimensional polyhedron of $\overline{\mathcal{M}}_{g, \Delta}$. The top-dimensional polyhedra exactly correspond to
combinatorial types which are $3$-valent and for which all marked points are in the interior of edges.
The image of the union of all top-dimensional polyhedra under the evaluation map does not yield a connected subset of $(\RR^2)^n$, but rather it divides the space of point conditions into cells which are separated from each other by walls. These walls are the images under the evaluation map of polyhedra of codimension $1$ in the moduli space. 

\begin{figure}
	\begin{center}
		\tikzset{every picture/.style={line width=0.75pt}} %set default line width to 0.75pt        

\begin{tikzpicture}[x=0.75pt,y=0.75pt,yscale=-1,xscale=1]
	%uncomment if require: \path (0,784); %set diagram left start at 0, and has height of 784
	
	%Straight Lines [id:da9012169808593964] 
	\draw    (200,140) -- (220,120) ;
	%Straight Lines [id:da4765703122477082] 
	\draw    (200,140) -- (220,160) ;
	%Straight Lines [id:da5526706469814957] 
	\draw    (150,140) -- (130,160) ;
	%Straight Lines [id:da6356030163975532] 
	\draw    (410,140) -- (370,140) ;
	%Straight Lines [id:da2470312105505632] 
	\draw  [dash pattern={on 0.84pt off 2.51pt}]  (120,170) -- (130,160) ;
	%Straight Lines [id:da37533549471217575] 
	\draw  [dash pattern={on 0.84pt off 2.51pt}]  (130.2,120) -- (130.2,110) ;
	%Straight Lines [id:da3679780668389995] 
	\draw  [dash pattern={on 0.84pt off 2.51pt}]  (120,120) -- (130,120) ;
	%Straight Lines [id:da006839717923191602] 
	\draw    (130.2,120) -- (150,140) ;
	%Curve Lines [id:da226249781128388] 
	\draw    (150,140) .. controls (170.2,130.07) and (180.2,130.87) .. (200,140) ;
	%Curve Lines [id:da8492078859407293] 
	\draw    (150,140) .. controls (169,149.27) and (180.6,148.87) .. (200,140) ;
	%Straight Lines [id:da9855230528264964] 
	\draw    (240,140) -- (288,140) ;
	\draw [shift={(290,140)}, rotate = 180] [color={rgb, 255:red, 0; green, 0; blue, 0 }  ][line width=0.75]    (10.93,-3.29) .. controls (6.95,-1.4) and (3.31,-0.3) .. (0,0) .. controls (3.31,0.3) and (6.95,1.4) .. (10.93,3.29)   ;
	%Straight Lines [id:da44914937668238253] 
	\draw    (410,140) -- (430,120) ;
	%Straight Lines [id:da7191646388225477] 
	\draw    (350,160) -- (370,140) ;
	%Straight Lines [id:da29525452500922356] 
	\draw    (430,160) -- (410,140) ;
	%Straight Lines [id:da5638150218186172] 
	\draw    (370,140) -- (350,120) ;
	%Straight Lines [id:da754700600290774] 
	\draw  [dash pattern={on 0.84pt off 2.51pt}]  (220,120) -- (230,110) ;
	%Straight Lines [id:da37680912333243066] 
	\draw  [dash pattern={on 0.84pt off 2.51pt}]  (220,160) -- (230,150) ;
	%Straight Lines [id:da7083364366305368] 
	\draw  [dash pattern={on 0.84pt off 2.51pt}]  (220,170) -- (220,160) ;
	%Straight Lines [id:da9117679974429109] 
	\draw  [dash pattern={on 0.84pt off 2.51pt}]  (430,170) -- (430,160) ;
	%Straight Lines [id:da7368700297222094] 
	\draw  [dash pattern={on 0.84pt off 2.51pt}]  (440,150) -- (430,160) ;
	%Straight Lines [id:da8372996876120572] 
	\draw  [dash pattern={on 0.84pt off 2.51pt}]  (430,120) -- (440,110) ;
	%Straight Lines [id:da6312431822149164] 
	\draw  [dash pattern={on 0.84pt off 2.51pt}]  (340,170) -- (350,160) ;
	%Straight Lines [id:da9144193462256865] 
	\draw  [dash pattern={on 0.84pt off 2.51pt}]  (350,120) -- (350,110) ;
	%Straight Lines [id:da32466627102480505] 
	\draw  [dash pattern={on 0.84pt off 2.51pt}]  (340,120) -- (350,120) ;
	%Shape: Rectangle [id:dp9138922564393571] 
	\draw   (310,90) -- (470,90) -- (470,200) -- (310,200) -- cycle ;
	
	% Text Node
	\draw (441,178) node [anchor=north west][inner sep=0.75pt]   [align=left] {$\displaystyle \mathbb{R}^{2}$};
	% Text Node
	\draw (182.2,125) node [anchor=north west][inner sep=0.75pt]   [align=left] {};
	% Text Node
	\draw (182.2,139) node [anchor=north west][inner sep=0.75pt]   [align=left] {};
	% Text Node
	\draw (269.5,127) node   [align=left] {\begin{minipage}[lt]{13.46pt}\setlength\topsep{0pt}
			$\displaystyle \varphi $
	\end{minipage}};
	% Text Node
	\draw (171,170) node [anchor=north west][inner sep=0.75pt]   [align=left] {$\displaystyle \Gamma $};

\end{tikzpicture}
	\end{center}
	\caption{A local picture of a parametrized tropical curve of exceptional type as in (3).}
	\label{fig-exceptional}
\end{figure}
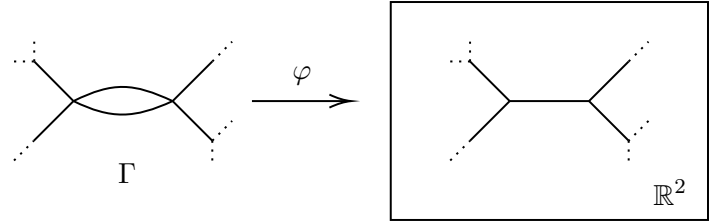

To show that the function counting preimages of the evaluation map with multiplicity is globally constant, we only need to cross these walls and show the local invariance around any wall. The walls have been classified in 
Theorem 4.8 \cite{GM051} and are given as 
%
%classifies these walls and shows the local invariance of the complex multiplicity for each type of wall. Section 4 in \cite{IKS09} shows the local invariance of the real multiplicity for each type of wall.
%
%Since the combinatorics of the tropical plane curve count stays the same no matter whether we are working on a complex count or an arithmetic count, we can prove Theorem \ref{thm-invariance} by showing the local invariance of the arithmetic multiplicity for each type of wall.
%Since the arithmetic count can be viewed as an interpolation between real and complex count, the proof uses both the real local invariance and the complex local invariance as ingredients.
%
%In order to show the arithmetic local invariance, let us first discuss the different types of walls. They are 
the images of the evaluation map of polyhedra corresponding to combinatorial types which
\begin{enumerate}
\item have precisely one $4$-valent vertex, all others $3$-valent, and all marked points in the interior of edges, or
\item are $3$ -valent and have precisely one marked point on a vertex and all others in the interior of edges, or
\item contain two edges that share the same two $4$-valent end vertices.
\end{enumerate}
The latter case is called an \emph{exceptional} combinatorial type, because the codimension of its associated polyhedron in the moduli space of tropical curves is not as expected by the formula given in Definition 3.3 of \cite{GM051}.
The two edges with the same end vertices must have the same primitive direction and are mapped to the same line segment in $\RR^2$, see Figure \ref{fig-exceptional}.

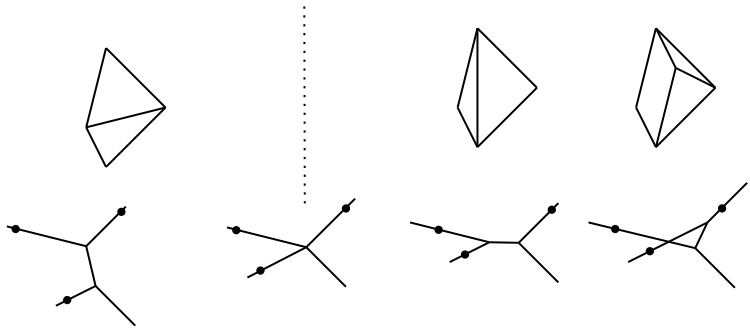
\begin{figure}
	\begin{center}
		\tikzset{every picture/.style={line width=0.75pt}} %set default line width to 0.75pt        

\begin{tikzpicture}[x=0.75pt,y=0.75pt,yscale=-1,xscale=1]
	%uncomment if require: \path (0,784); %set diagram left start at 0, and has height of 784
	
	%Straight Lines [id:da9012169808593964] 
	\draw    (150,310) -- (140,290) ;
	%Straight Lines [id:da4765703122477082] 
	\draw    (150,250) -- (180,280) ;
	%Straight Lines [id:da5526706469814957] 
	\draw    (150,250) -- (140,290) ;
	%Straight Lines [id:da006839717923191602] 
	\draw    (180,280) -- (150,310) ;
	%Straight Lines [id:da42090316851078036] 
	\draw    (427.33,300) -- (417.33,280) ;
	%Straight Lines [id:da7417224989885193] 
	\draw    (427.33,240) -- (457.33,270) ;
	%Straight Lines [id:da13696837670762252] 
	\draw    (427.33,240) -- (417.33,280) ;
	%Straight Lines [id:da800121658840217] 
	\draw    (457.33,270) -- (427.33,300) ;
	%Straight Lines [id:da3505592333715972] 
	\draw    (337.33,300) -- (327.33,280) ;
	%Straight Lines [id:da48101539029986673] 
	\draw    (337.33,240) -- (367.33,270) ;
	%Straight Lines [id:da086064008776221] 
	\draw    (337.33,240) -- (327.33,280) ;
	%Straight Lines [id:da2837118019182576] 
	\draw    (367.33,270) -- (337.33,300) ;
	%Straight Lines [id:da7473463236750383] 
	\draw  [dash pattern={on 0.84pt off 2.51pt}]  (250,229) -- (250.25,329.72) ;
	%Straight Lines [id:da3407586093459145] 
	\draw    (437.33,260) -- (427.33,240) ;
	%Straight Lines [id:da9152090551192767] 
	\draw    (437.33,260) -- (427.33,300) ;
	%Straight Lines [id:da8613333967440306] 
	\draw    (457.33,270) -- (437.33,260) ;
	%Straight Lines [id:da3698121622321914] 
	\draw    (337.33,300) -- (337.33,240) ;
	%Straight Lines [id:da8708579289334668] 
	\draw    (140,290) -- (180,280) ;
	%Straight Lines [id:da6719490959786504] 
	\draw    (140,350) -- (100,340) ;
	%Straight Lines [id:da0014040713410351513] 
	\draw    (160,330) -- (140,350) ;
	%Straight Lines [id:da9256307362175225] 
	\draw    (140,350) -- (144.73,370.07) ;
	%Straight Lines [id:da7624303518639748] 
	\draw    (124.73,380.07) -- (144.73,370.07) ;
	%Straight Lines [id:da8584529767983721] 
	\draw    (144.73,370.07) -- (164.73,390.07) ;
	%Straight Lines [id:da11833551589916302] 
	\draw    (343.33,348) -- (303.33,338) ;
	%Straight Lines [id:da5715447102001585] 
	\draw    (447.17,350.97) -- (393.33,338) ;
	%Straight Lines [id:da09419526182418325] 
	\draw    (378.13,328.27) -- (358.13,348.27) ;
	%Straight Lines [id:da7067148637852558] 
	\draw    (473.33,318) -- (453.33,338) ;
	%Straight Lines [id:da30383989892904917] 
	\draw    (323.33,358) -- (343.33,348) ;
	%Straight Lines [id:da541662809668358] 
	\draw    (413.33,358) -- (453.33,338) ;
	%Straight Lines [id:da4435546969112272] 
	\draw    (358.13,348.27) -- (378.13,368.27) ;
	%Straight Lines [id:da01383102518387791] 
	\draw    (447.17,350.97) -- (467.17,370.97) ;
	%Straight Lines [id:da7870333723838169] 
	\draw    (343.33,348) -- (358.13,348.27) ;
	%Straight Lines [id:da11074234195346644] 
	\draw    (453.33,338) -- (447.17,350.97) ;
	%Straight Lines [id:da3990248721735148] 
	\draw    (251,350.5) -- (211,340.5) ;
	%Straight Lines [id:da3891735828760884] 
	\draw    (275.62,325.97) -- (251,350.5) ;
	%Straight Lines [id:da6864554498913713] 
	\draw    (251,350.5) -- (271,370.5) ;
	%Straight Lines [id:da1386372812614961] 
	\draw    (221.31,365.81) -- (251,350.5) ;
	%Shape: Circle [id:dp0843269101581392] 
	\draw  [fill={rgb, 255:red, 0; green, 0; blue, 0 }  ,fill opacity=1 ] (316.17,341.72) .. controls (316.17,340.84) and (316.88,340.13) .. (317.75,340.13) .. controls (318.62,340.13) and (319.33,340.84) .. (319.33,341.72) .. controls (319.33,342.59) and (318.62,343.3) .. (317.75,343.3) .. controls (316.88,343.3) and (316.17,342.59) .. (316.17,341.72) -- cycle ;
	%Shape: Circle [id:dp8407287032972187] 
	\draw  [fill={rgb, 255:red, 0; green, 0; blue, 0 }  ,fill opacity=1 ] (329.48,354.18) .. controls (329.48,353.31) and (330.19,352.6) .. (331.07,352.6) .. controls (331.94,352.6) and (332.65,353.31) .. (332.65,354.18) .. controls (332.65,355.06) and (331.94,355.77) .. (331.07,355.77) .. controls (330.19,355.77) and (329.48,355.06) .. (329.48,354.18) -- cycle ;
	%Shape: Circle [id:dp8939996788434548] 
	\draw  [fill={rgb, 255:red, 0; green, 0; blue, 0 }  ,fill opacity=1 ] (405.23,341.32) .. controls (405.23,340.44) and (405.94,339.73) .. (406.82,339.73) .. controls (407.69,339.73) and (408.4,340.44) .. (408.4,341.32) .. controls (408.4,342.19) and (407.69,342.9) .. (406.82,342.9) .. controls (405.94,342.9) and (405.23,342.19) .. (405.23,341.32) -- cycle ;
	%Shape: Circle [id:dp9984092454025638] 
	\draw  [fill={rgb, 255:red, 0; green, 0; blue, 0 }  ,fill opacity=1 ] (422.7,352.52) .. controls (422.7,351.64) and (423.41,350.93) .. (424.28,350.93) .. controls (425.16,350.93) and (425.87,351.64) .. (425.87,352.52) .. controls (425.87,353.39) and (425.16,354.1) .. (424.28,354.1) .. controls (423.41,354.1) and (422.7,353.39) .. (422.7,352.52) -- cycle ;
	%Shape: Circle [id:dp4938912157210078] 
	\draw  [fill={rgb, 255:red, 0; green, 0; blue, 0 }  ,fill opacity=1 ] (459.1,330.92) .. controls (459.1,330.04) and (459.81,329.33) .. (460.68,329.33) .. controls (461.56,329.33) and (462.27,330.04) .. (462.27,330.92) .. controls (462.27,331.79) and (461.56,332.5) .. (460.68,332.5) .. controls (459.81,332.5) and (459.1,331.79) .. (459.1,330.92) -- cycle ;
	%Shape: Circle [id:dp18408884082211185] 
	\draw  [fill={rgb, 255:red, 0; green, 0; blue, 0 }  ,fill opacity=1 ] (373.23,331.58) .. controls (373.23,330.71) and (373.94,330) .. (374.82,330) .. controls (375.69,330) and (376.4,330.71) .. (376.4,331.58) .. controls (376.4,332.46) and (375.69,333.17) .. (374.82,333.17) .. controls (373.94,333.17) and (373.23,332.46) .. (373.23,331.58) -- cycle ;
	%Shape: Circle [id:dp03447555849288342] 
	\draw  [fill={rgb, 255:red, 0; green, 0; blue, 0 }  ,fill opacity=1 ] (269.44,330.93) .. controls (269.44,330.06) and (270.15,329.35) .. (271.03,329.35) .. controls (271.9,329.35) and (272.61,330.06) .. (272.61,330.93) .. controls (272.61,331.81) and (271.9,332.52) .. (271.03,332.52) .. controls (270.15,332.52) and (269.44,331.81) .. (269.44,330.93) -- cycle ;
	%Shape: Circle [id:dp8606007917088517] 
	\draw  [fill={rgb, 255:red, 0; green, 0; blue, 0 }  ,fill opacity=1 ] (214.11,341.9) .. controls (214.11,341.02) and (214.82,340.32) .. (215.69,340.32) .. controls (216.57,340.32) and (217.28,341.02) .. (217.28,341.9) .. controls (217.28,342.77) and (216.57,343.48) .. (215.69,343.48) .. controls (214.82,343.48) and (214.11,342.77) .. (214.11,341.9) -- cycle ;
	%Shape: Circle [id:dp3562983981370962] 
	\draw  [fill={rgb, 255:red, 0; green, 0; blue, 0 }  ,fill opacity=1 ] (226.31,362.3) .. controls (226.31,361.42) and (227.01,360.72) .. (227.89,360.72) .. controls (228.76,360.72) and (229.47,361.42) .. (229.47,362.3) .. controls (229.47,363.17) and (228.76,363.88) .. (227.89,363.88) .. controls (227.01,363.88) and (226.31,363.17) .. (226.31,362.3) -- cycle ;
	%Shape: Circle [id:dp10203872702140482] 
	\draw  [fill={rgb, 255:red, 0; green, 0; blue, 0 }  ,fill opacity=1 ] (102.83,341.32) .. controls (102.83,340.44) and (103.54,339.73) .. (104.42,339.73) .. controls (105.29,339.73) and (106,340.44) .. (106,341.32) .. controls (106,342.19) and (105.29,342.9) .. (104.42,342.9) .. controls (103.54,342.9) and (102.83,342.19) .. (102.83,341.32) -- cycle ;
	%Shape: Circle [id:dp10079406868909813] 
	\draw  [fill={rgb, 255:red, 0; green, 0; blue, 0 }  ,fill opacity=1 ] (156.17,332.65) .. controls (156.17,331.78) and (156.88,331.07) .. (157.75,331.07) .. controls (158.62,331.07) and (159.33,331.78) .. (159.33,332.65) .. controls (159.33,333.52) and (158.62,334.23) .. (157.75,334.23) .. controls (156.88,334.23) and (156.17,333.52) .. (156.17,332.65) -- cycle ;
	%Shape: Circle [id:dp22197016930116054] 
	\draw  [fill={rgb, 255:red, 0; green, 0; blue, 0 }  ,fill opacity=1 ] (128.83,377.18) .. controls (128.83,376.31) and (129.54,375.6) .. (130.42,375.6) .. controls (131.29,375.6) and (132,376.31) .. (132,377.18) .. controls (132,378.06) and (131.29,378.77) .. (130.42,378.77) .. controls (129.54,378.77) and (128.83,378.06) .. (128.83,377.18) -- cycle ;

\end{tikzpicture}
	\end{center}
	\caption{The picture shows a wall of type (1) together with local pictures of the tropical curves passing through the points on the wall, as well as on both sides of the walls, and their dual Newton subdivisions. The depicted points should be viewed as an example of how the point conditions and balancing could fix the curves. In this example, the wall crossing happens if we move the lowest point upwards.}
	\label{fig-4valent}
\end{figure}

For case (1), Figure \ref{fig-4valent} shows a local picture of the three resolutions of a tropical curve with one $4$-valent vertex. As all other vertices and their multiplicity stay the same in all three resolutions, by Lemma~\ref{lem-vertexmult} it is enough to study this local situation and compare the products of the multiplicities of the two vertices in each. One of the three resolutions produces a crossing of two edges corresponding to a paralellogram in its dual Newton subdivision. It follows from Theorem 4.8 \cite{GM051} and Section 4.6 \cite{IKS09} that the tropical curve dual to the subdivision with the parallelogram occurs on the same side of the wall as the one containing the diagonal of the parallelogram, see Figure \ref{fig-4valent}.

\begin{figure}
	\begin{center}
		\tikzset{every picture/.style={line width=0.75pt}} %set default line width to 0.75pt        

\begin{tikzpicture}[x=0.75pt,y=0.75pt,yscale=-1,xscale=1]
	%uncomment if require: \path (0,784); %set diagram left start at 0, and has height of 784
	
	%Straight Lines [id:da7473463236750383] 
	\draw  [dash pattern={on 0.84pt off 2.51pt}]  (250,229) -- (250.25,329.72) ;
	%Straight Lines [id:da6719490959786504] 
	\draw    (190,280) -- (150,280) ;
	%Straight Lines [id:da0014040713410351513] 
	\draw    (220,250) -- (190,280) ;
	%Straight Lines [id:da9256307362175225] 
	\draw    (190,280) -- (190,320) ;
	%Straight Lines [id:da3891735828760884] 
	\draw    (340,250) -- (310,280) ;
	%Shape: Circle [id:dp0843269101581392] 
	\draw  [fill={rgb, 255:red, 0; green, 0; blue, 0 }  ,fill opacity=1 ] (308.42,301.58) .. controls (308.42,300.71) and (309.13,300) .. (310,300) .. controls (310.87,300) and (311.58,300.71) .. (311.58,301.58) .. controls (311.58,302.46) and (310.87,303.17) .. (310,303.17) .. controls (309.13,303.17) and (308.42,302.46) .. (308.42,301.58) -- cycle ;
	%Shape: Circle [id:dp18408884082211185] 
	\draw  [fill={rgb, 255:red, 0; green, 0; blue, 0 }  ,fill opacity=1 ] (188.43,302.98) .. controls (188.43,302.11) and (189.14,301.4) .. (190.02,301.4) .. controls (190.89,301.4) and (191.6,302.11) .. (191.6,302.98) .. controls (191.6,303.86) and (190.89,304.57) .. (190.02,304.57) .. controls (189.14,304.57) and (188.43,303.86) .. (188.43,302.98) -- cycle ;
	%Shape: Circle [id:dp03447555849288342] 
	\draw  [fill={rgb, 255:red, 0; green, 0; blue, 0 }  ,fill opacity=1 ] (321.84,266.53) .. controls (321.84,265.66) and (322.55,264.95) .. (323.43,264.95) .. controls (324.3,264.95) and (325.01,265.66) .. (325.01,266.53) .. controls (325.01,267.41) and (324.3,268.12) .. (323.43,268.12) .. controls (322.55,268.12) and (321.84,267.41) .. (321.84,266.53) -- cycle ;
	%Shape: Circle [id:dp10079406868909813] 
	\draw  [fill={rgb, 255:red, 0; green, 0; blue, 0 }  ,fill opacity=1 ] (170,280) .. controls (170,279.13) and (170.71,278.42) .. (171.58,278.42) .. controls (172.46,278.42) and (173.17,279.13) .. (173.17,280) .. controls (173.17,280.87) and (172.46,281.58) .. (171.58,281.58) .. controls (170.71,281.58) and (170,280.87) .. (170,280) -- cycle ;
	%Straight Lines [id:da054606158417235484] 
	\draw    (310,280) -- (310,320) ;
	%Straight Lines [id:da11663753710733504] 
	\draw    (310,280) -- (270,280) ;
	%Straight Lines [id:da9312911290559469] 
	\draw    (250,370) -- (210,370) ;
	%Straight Lines [id:da2494590579528475] 
	\draw    (280,340) -- (250,370) ;
	%Straight Lines [id:da05884037160143607] 
	\draw    (250,370) -- (250,410) ;
	%Shape: Circle [id:dp31342397051602655] 
	\draw  [fill={rgb, 255:red, 0; green, 0; blue, 0 }  ,fill opacity=1 ] (248.42,391.58) .. controls (248.42,390.71) and (249.13,390) .. (250,390) .. controls (250.87,390) and (251.58,390.71) .. (251.58,391.58) .. controls (251.58,392.46) and (250.87,393.17) .. (250,393.17) .. controls (249.13,393.17) and (248.42,392.46) .. (248.42,391.58) -- cycle ;
	%Shape: Circle [id:dp2471587919511412] 
	\draw  [fill={rgb, 255:red, 0; green, 0; blue, 0 }  ,fill opacity=1 ] (248.42,370) .. controls (248.42,369.13) and (249.13,368.42) .. (250,368.42) .. controls (250.87,368.42) and (251.58,369.13) .. (251.58,370) .. controls (251.58,370.87) and (250.87,371.58) .. (250,371.58) .. controls (249.13,371.58) and (248.42,370.87) .. (248.42,370) -- cycle ;

\end{tikzpicture}
	\end{center}
	\caption{The picture shows a wall of type (2) together with local pictures of the tropical curves passing through the points. On the wall, there is a point on a vertex. The other point should be viewed only as an example of how the remaining points fix the curve. For this example, we cross the wall by moving the top point to the right.}
	\label{fig-pointonvertex}
\end{figure}
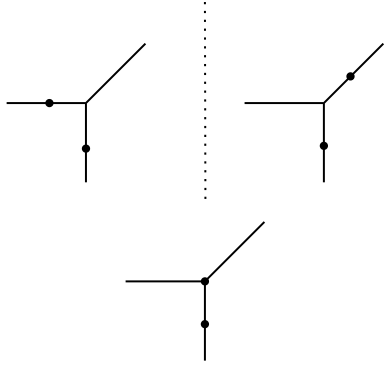

For an example of case (2) see Figure~\ref{fig-pointonvertex}. Let $(\Gamma,\varphi)$ be a parametrized tropical curve through points $p_1,\ldots,p_n$ in general position, where $n=\#\Delta+g-1$. Then by Lemma~4.20 \cite{Mi03} we have that $\Gamma$ minus the inverse images of the $p_i$ under $\varphi$ is a forest such that each connected component contains precisely one end. Let us assume that $p_1$ is the point on a vertex, and consider a tropical curve $(\Gamma,\varphi)$ which lives on one side of the wall. Its combinatorial type corresponds to a polyhedron of top-dimension in the moduli space, i.e.\ $p_1$ must have wandered inside an edge. We consider $\Gamma\setminus\big( \varphi^{-1}(p_2)\cup\ldots,\cup \varphi^{-1}(p_n)\big)$. By the above, it has one connected component which is either a tree but contains two ends, or one connected component of genus $1$. This connected component contains precisely two of the edges adjacent to the vertex containing $p_1$ after specializing the point conditions to live on our wall. Accordingly, the point can only wander on two of the adjacent edges and not on the third, see Figure \ref{fig-pointonvertex}. It is clear that those two resolutions live on opposite sides of the wall.

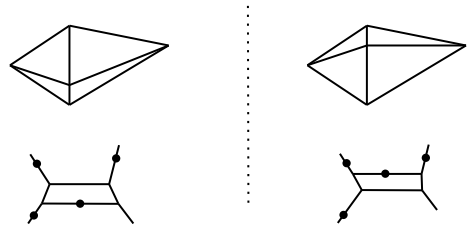
\begin{figure}
	\begin{center}
		\tikzset{every picture/.style={line width=0.75pt}} %set default line width to 0.75pt        

\begin{tikzpicture}[x=0.75pt,y=0.75pt,yscale=-1,xscale=1]
	%uncomment if require: \path (0,784); %set diagram left start at 0, and has height of 784
	
	%Straight Lines [id:da7473463236750383] 
	\draw  [dash pattern={on 0.84pt off 2.51pt}]  (250,140) -- (250.25,240.72) ;
	%Straight Lines [id:da6719490959786504] 
	\draw    (160,150) -- (130,170) ;
	%Straight Lines [id:da0014040713410351513] 
	\draw    (210,160) -- (160,150) ;
	%Straight Lines [id:da9256307362175225] 
	\draw    (160,150) -- (160,190) ;
	%Straight Lines [id:da3891735828760884] 
	\draw    (185,210.31) -- (180,230) ;
	%Shape: Circle [id:dp18408884082211185] 
	\draw  [fill={rgb, 255:red, 0; green, 0; blue, 0 }  ,fill opacity=1 ] (181.97,216.98) .. controls (181.97,216.11) and (182.68,215.4) .. (183.56,215.4) .. controls (184.43,215.4) and (185.14,216.11) .. (185.14,216.98) .. controls (185.14,217.86) and (184.43,218.57) .. (183.56,218.57) .. controls (182.68,218.57) and (181.97,217.86) .. (181.97,216.98) -- cycle ;
	%Shape: Circle [id:dp10079406868909813] 
	\draw  [fill={rgb, 255:red, 0; green, 0; blue, 0 }  ,fill opacity=1 ] (142,219.69) .. controls (142,218.82) and (142.71,218.11) .. (143.58,218.11) .. controls (144.46,218.11) and (145.17,218.82) .. (145.17,219.69) .. controls (145.17,220.57) and (144.46,221.28) .. (143.58,221.28) .. controls (142.71,221.28) and (142,220.57) .. (142,219.69) -- cycle ;
	%Straight Lines [id:da054606158417235484] 
	\draw    (150,230) -- (180,230) ;
	%Straight Lines [id:da11663753710733504] 
	\draw    (150,230) -- (140.23,214.93) ;
	%Straight Lines [id:da9312911290559469] 
	\draw    (150,230) -- (146.08,239.77) ;
	%Straight Lines [id:da2494590579528475] 
	\draw    (180,230) -- (184.69,239.93) ;
	%Straight Lines [id:da05884037160143607] 
	\draw    (146.08,239.77) -- (184.69,239.93) ;
	%Shape: Circle [id:dp31342397051602655] 
	\draw  [fill={rgb, 255:red, 0; green, 0; blue, 0 }  ,fill opacity=1 ] (163.8,239.85) .. controls (163.8,238.98) and (164.51,238.27) .. (165.38,238.27) .. controls (166.26,238.27) and (166.97,238.98) .. (166.97,239.85) .. controls (166.97,240.73) and (166.26,241.43) .. (165.38,241.43) .. controls (164.51,241.43) and (163.8,240.73) .. (163.8,239.85) -- cycle ;
	%Shape: Circle [id:dp2471587919511412] 
	\draw  [fill={rgb, 255:red, 0; green, 0; blue, 0 }  ,fill opacity=1 ] (140.46,245.8) .. controls (140.46,244.93) and (141.17,244.22) .. (142.04,244.22) .. controls (142.92,244.22) and (143.63,244.93) .. (143.63,245.8) .. controls (143.63,246.68) and (142.92,247.38) .. (142.04,247.38) .. controls (141.17,247.38) and (140.46,246.68) .. (140.46,245.8) -- cycle ;
	%Straight Lines [id:da4457526209314473] 
	\draw    (130,170) -- (160,190) ;
	%Straight Lines [id:da3890292067693397] 
	\draw    (210,160) -- (160,190) ;
	%Straight Lines [id:da031407231532135205] 
	\draw    (130,170) -- (160,180) ;
	%Straight Lines [id:da9825597184491929] 
	\draw    (160,180) -- (210,160) ;
	%Straight Lines [id:da5794101547140469] 
	\draw    (310,150) -- (280,170) ;
	%Straight Lines [id:da5874802066682225] 
	\draw    (360,160) -- (310,150) ;
	%Straight Lines [id:da4682907003649821] 
	\draw    (310,150) -- (310,190) ;
	%Straight Lines [id:da08617904874035742] 
	\draw    (280,170) -- (310,190) ;
	%Straight Lines [id:da49331758820371974] 
	\draw    (360,160) -- (310,190) ;
	%Straight Lines [id:da17045475516490005] 
	\draw    (280,170) -- (310,160) ;
	%Straight Lines [id:da3056691784946226] 
	\draw    (310,160) -- (360,160) ;
	%Straight Lines [id:da8422418106503133] 
	\draw    (139.15,249.85) -- (146.08,239.77) ;
	%Straight Lines [id:da04578395722550466] 
	\draw    (192.23,249.85) -- (184.69,239.93) ;
	%Straight Lines [id:da5264366752518663] 
	\draw    (341.18,210.03) -- (337.57,224.82) ;
	%Shape: Circle [id:dp5229807227855501] 
	\draw  [fill={rgb, 255:red, 0; green, 0; blue, 0 }  ,fill opacity=1 ] (338.15,216.7) .. controls (338.15,215.83) and (338.86,215.12) .. (339.73,215.12) .. controls (340.61,215.12) and (341.32,215.83) .. (341.32,216.7) .. controls (341.32,217.58) and (340.61,218.28) .. (339.73,218.28) .. controls (338.86,218.28) and (338.15,217.58) .. (338.15,216.7) -- cycle ;
	%Shape: Circle [id:dp16045602296407369] 
	\draw  [fill={rgb, 255:red, 0; green, 0; blue, 0 }  ,fill opacity=1 ] (298.18,219.41) .. controls (298.18,218.54) and (298.89,217.83) .. (299.76,217.83) .. controls (300.64,217.83) and (301.35,218.54) .. (301.35,219.41) .. controls (301.35,220.28) and (300.64,220.99) .. (299.76,220.99) .. controls (298.89,220.99) and (298.18,220.28) .. (298.18,219.41) -- cycle ;
	%Straight Lines [id:da8734438454937499] 
	\draw    (302.93,224.84) -- (337.57,224.82) ;
	%Straight Lines [id:da04106830613829693] 
	\draw    (302.93,224.84) -- (296.41,214.65) ;
	%Straight Lines [id:da02102717769015483] 
	\draw    (302.93,224.84) -- (307.43,232.97) ;
	%Straight Lines [id:da020707151120656686] 
	\draw    (337.57,224.82) -- (337.86,233.11) ;
	%Straight Lines [id:da24650318004109484] 
	\draw    (307.43,232.97) -- (337.86,233.11) ;
	%Shape: Circle [id:dp06001031495372244] 
	\draw  [fill={rgb, 255:red, 0; green, 0; blue, 0 }  ,fill opacity=1 ] (317.73,224.72) .. controls (317.73,223.84) and (318.44,223.13) .. (319.31,223.13) .. controls (320.18,223.13) and (320.89,223.84) .. (320.89,224.72) .. controls (320.89,225.59) and (320.18,226.3) .. (319.31,226.3) .. controls (318.44,226.3) and (317.73,225.59) .. (317.73,224.72) -- cycle ;
	%Shape: Circle [id:dp3856850866046483] 
	\draw  [fill={rgb, 255:red, 0; green, 0; blue, 0 }  ,fill opacity=1 ] (296.64,245.52) .. controls (296.64,244.64) and (297.35,243.94) .. (298.22,243.94) .. controls (299.1,243.94) and (299.81,244.64) .. (299.81,245.52) .. controls (299.81,246.39) and (299.1,247.1) .. (298.22,247.1) .. controls (297.35,247.1) and (296.64,246.39) .. (296.64,245.52) -- cycle ;
	%Straight Lines [id:da32697040972387525] 
	\draw    (295.33,249.57) -- (307.43,232.97) ;
	%Straight Lines [id:da9421113081235678] 
	\draw    (345.4,243.03) -- (337.86,233.11) ;

\end{tikzpicture}
	\end{center}
	\caption{The picture shows a wall of type (3) together with local pictures of the tropical curves passing through points on both sides of the wall. Again, the points should be viewed as an example. For this example, we cross the wall by moving the point on the cycle upwards.}
	\label{fig-resolveexceptional}
\end{figure}

For case (3) finally, Figure \ref{fig-resolveexceptional} shows the two possible resolutions. Note that even if the two parallel edges have the same weight, they are distinguishable as at least one edge must contain a point, which follows from the generic position of the points away from the wall and again Lemma 4.20 \cite{Mi03}.

\begin{proof}[Proof of Theorem \ref{thm-invariance}]
We show the local invariance of the arithmetic multiplicity for the three cases (1), (2), and (3) described above:
%\begin{enumerate}
%\item 

\textbf{Wall of type (1):} We assume that the $4$-valent vertex is non-degenerate, i.e.\ all three resolutions exist as sketched in Figure \ref{fig-4valent} and none contains an edge of direction $0$ or vertices whose adjacent edges are all parallel. All of these degenerate cases are similar, but potentially with fewer terms. 

Denote the resolution on the left of Figure~\ref{fig-4valent} by $C_L$ and the two resolutions on the right by $C_{R_1}$ and $C_{R_2}$. We have to show
\begin{equation} \label{eq-toshow}
	\mult_{\AA^1}(C_L) = \mult_{\AA^1}(C_{R_1}) + \mult_{\AA^1}(C_{R_2}).
\end{equation}
We present this calculation only in case that none of the four edges adjacent to the $4$-valent vertex is even, and that not all three newly emerging bounded edges are even. Any other case is completely hyperbolic, which means that the computation is similar but simpler.

Under our assumptions, it follows from Lemma~19 of \cite{IKS09} that precisely one of the newly emerging bounded edges is even while the other two are odd.
Assume first that the diagonal of the quadrilateral dual to the $4$-valent vertex which we inserted for the left resolution in Figure \ref{fig-4valent} is even. 
In this case \eqref{eq-toshow} becomes
\begin{equation} \label{eq-wall1}
\begin{aligned} 
	&\frac{\mult_{\mathbb{C}}(C_L)}{2}\cdot \mathbb{H}\\ 
	={}& \frac{\mult_{\mathbb{C}}(C_{R_1})-1}{2}\cdot \mathbb{H}+ \Big\langle (-1)^{i_1} \cdot \prod_{j=1}^4 w_j \Big\rangle  
	+ \frac{\mult_{\mathbb{C}}(C_{R_2})-1}{2}\cdot \mathbb{H}+ \Big\langle (-1)^{i_2} \cdot \prod_{j=1}^4 w_j \Big\rangle ,
\end{aligned}
\end{equation}
where $w_1,\ldots,w_4$ denote the weights of the $4$ exterior edges, and $i_1$ resp.\ $i_2$ denotes the number of interior points in triangles in the Newton subdivision dual to $C_{R_1}$ resp.\ $C_{R_2}$.
As the real count is invariant by Section 4 \cite{IKS09}, and we have an even edge on the left yielding a contribution of $0$ there, 
one of $(-1)^{i_1}$ and $(-1)^{i_2}$ must be negative, the other positive. Therefore, by Lemma~\ref{lem:hyperbolic} we have
\[ \Big\langle (-1)^{i_1} \cdot \prod_{j=1}^4 w_j \Big\rangle + \Big\langle (-1)^{i_2} \cdot \prod_{j=1}^4 w_j \Big\rangle = \mathbb{H}. \]
The equality in~\eqref{eq-wall1} then follows because of the invariance of the complex count (Theorem 4.8 \cite{GM051}), i.e. $\mult_{\mathbb{C}}(C_{R_1})+\mult_{\mathbb{C}}(C_{R_2})= \mult_{\mathbb{C}}(C_L)$.

Assume next that the diagonal of the quadrilateral in the middle resolution in Figure \ref{fig-4valent} is even.
We then have to show
\begin{align*} 
	&\frac{\mult_{\mathbb{C}}(C_L)-1}{2}\cdot \mathbb{H}+ \Big\langle (-1)^{i_0} \cdot \prod_{j=1}^4 w_j \Big\rangle \\ 
	={}& \frac{\mult_{\mathbb{C}}(C_{R_1})}{2}\cdot \mathbb{H}  + \frac{\mult_{\mathbb{C}}(C_{R_2})-1}{2}\cdot \mathbb{H}+ \Big\langle (-1)^{i_2} \cdot \prod_{j=1}^4 w_j \Big\rangle ,
\end{align*}
where $i_0$ denotes the number of interior points in triangles in the Newton subdivision dual to $C_L$.
By the real invariance, $(-1)^{i_0}$ must equal $(-1)^{i_2}$ in this case. The equality then follows from the complex invariance. 
For the final case, we use the same argument with the roles of $C_{R_1}$ and $C_{R_2}$ exchanged.

\textbf{Wall of type (2):} This case is clear, because the arithmetic multiplicity of a tropical curve only depends on the dual Newton subdivision of the curve, which is the same on both sides of a wall of type~(2).

\textbf{Wall of type (3):} It follows from the real invariance that all edges on the  left are odd if and only if all edges on the right are odd, see Section 4.5 in \cite{IKS09}. Furthermore, the real invariance shows that the signs on both sides must agree. The coefficients of $\mathbb{H}$ agree because of the complex invariance.

\end{proof}

\section{Caporaso-Harris formula}

Throughout this section we study plane tropical curves with Newton polygon $\Delta_d = \conv \big\{(0,0), \allowbreak (0, d), (d, 0) \big\}$ for some degree $d$. This restriction is mainly to simplify notation and a more general version is possible with the same underlying ideas.
The goal of this section is to prove Theorem~\ref{thm-caporaso-harris} and thereby establish the Caporaso-Harris formula to recursively compute the invariants $N^\trop_{\mathbb{A}^1}(d,g)$. To this end, we first have to extend our arithmetic invariants to allow counts of tropical curves with left ends of higher weight and left ends which are fixed by point conditions on the far left.
We fix some notation.
\begin{definition}
  A (finite) sequence is a collection $ \alpha=(\alpha_1,\alpha_2,\dots)
  $ of natural numbers almost all of which are zero. If $ \alpha_k=0 $ for all
  $ k>n $ we will also write this sequence as $
  \alpha=(\alpha_1,\dots,\alpha_n) $. For two sequences $ \alpha $ and $ \beta
  $ we define
    \[ \begin {array}{r@{\;}l@{\;}l}
      |\alpha| &:=& \alpha_1+\alpha_2+\cdots; \\
      I\alpha &:=& 1\alpha_1+2\alpha_2+3\alpha_3+\cdots; \\
      I^{\alpha} &:=& 1^{\alpha_1}\cdot 2^{\alpha_2}\cdot 3^{\alpha_3}\cdot
        \; \cdots ; \\
      \alpha+\beta &:=& (\alpha_1+\beta_1,\alpha_2+\beta_2,\ldots); \\
      \alpha \geq \beta &:\Leftrightarrow& \alpha_n \geq \beta_n
        \mbox { for all $n$}; \\
      \binom{\alpha}{\beta} &:=& \binom{\alpha_1}{\beta_1}\cdot
        \binom{\alpha_2}{\beta_2}\cdot \; \cdots.
    \end {array} \]
  We denote by $e_k$ the sequence which has a $1$ at the $k$-th place and zeros
  everywhere else.
\end{definition}

\begin{definition} \label{def-arithmultwithfixedends}
  Let $C$ be a simple tropical curve of degree $d$ and genus $g$ with
  $\alpha_i$ fixed and $\beta_i$ non-fixed ends to the left of weight
  $i$ for all $i$. We define the \emph{arithmetic $(\alpha,\beta)$-multiplicity} of $C$ to
  be
    $$ \mult_{\mathbb{A}^1}^{\alpha,\beta}(C)
    =\begin{cases}
    	\frac{1}{2}\big(\frac{1}{I^{\alpha}}  \mult_{\mathbb{C}}(C)-1\big)\cdot \mathbb{H}
    	+\big\langle (-1)^i I^\beta \big\rangle & \mbox{if $\frac{1}{I^{\alpha}} \mult_{\mathbb{C}}(C)$ is odd,} \\ 
    	\frac{1}{2} \big(\frac{1}{I^{\alpha}} \mult_{\mathbb{C}}(C)\big)\cdot \mathbb{H}&\mbox{if $\frac{1}{I^{\alpha}} \mult_{\mathbb{C}}(C)$ is even,}
    \end{cases} $$
  where $i$ denotes the number of interior lattice points in triangles in the Newton subdivision dual to $C$.

  Let $ d, g \ge 0 $ be integers, and let $ \alpha $ and $ \beta $ be
  sequences with $ I\alpha+I\beta=d $. Then we define $N_{\mathbb{A}^1}^{\trop,\alpha,\beta}(d,g)$ to be the number of tropical curves of
  degree $d$ and genus $g$ with $\alpha_i$ fixed and $\beta_i$ non-fixed
   ends to the left of weight $i$ for all $i$ that pass in addition
  through a set of $ 2d+g+|\beta|-1 $ points in general position.
  The curves are to be counted with their respective arithmetic $(\alpha,\beta)
  $-multiplicities. 
\end{definition}

By Theorem \ref{thm-invariance}, the definition of $N_{\mathbb{A}^1}^{\trop,\alpha,\beta}(d,g)$ does not depend on the
  configuration of point conditions and fixed ends.
 Using the techniques involved in the proof of the Correspondence Theorem \ref{thm-corres-new}, the numbers $N_{\mathbb{A}^1}^{\trop,\alpha,\beta}(d,g)$ can be shown to equal an arithmetic count of algebraic curves of degree $d$ and genus $g$ close to the tropical limit, passing through the right number of points, and with tangency conditions to a fixed line encoded by the sequences $\alpha$ and $\beta$ as follows: there must be $\alpha_i$ points of contact order $i$ in fixed points on the line and $\beta_i$ further points of contact order $i$.

\begin{theorem} \label{thm-caporaso-harris}
Let $K$ be a field of characteristic larger than $d$.
The numbers $ N_{\mathbb{A}^1}^{\trop,\alpha,\beta}(d,g)\in \GW(K)$  satisfy the following Caporaso-Harris recursion:
  
  \begin{equation} \label{eq:Caporaso_Harris}
  \begin{aligned}
    N_{\mathbb{A}^1}^{\trop,\alpha,\beta}(d,g)
      =& \sum_{k \mbox{ \tiny odd}:\beta_k>0}\Big( \frac{k-1}{2}\cdot \mathbb{H}+\langle k \rangle\Big) \cdot N_{\mathbb{A}^1}^{\trop,\alpha+e_k,\beta-e_k}(d,g)
        \\
        &  +\sum_{k \mbox{ \tiny  even}:\beta_k>0} \Big( \frac{k}{2}\cdot \mathbb{H}\Big)\cdot N_{\mathbb{A}^1}^{\trop,\alpha+e_k,\beta-e_k}(d,g)
        \\
      &+ \sum\Big(\frac{ I^{\beta'-\beta}-1}{2}\cdot \mathbb{H}+ \langle I^{\beta'-\beta}\rangle\Big)  \cdot \binom{\alpha}{\alpha'} \cdot
        \binom{\beta'}{\beta} \cdot N_{\mathbb{A}^1}^{\trop,\alpha',\beta'}(d-1,g')\\
       &  + \sum\Big(\frac{ I^{\beta'-\beta}-1}{2}\cdot \mathbb{H}\Big)  \cdot \binom{\alpha}{\alpha'} \cdot
        \binom{\beta'}{\beta} \cdot N_{\mathbb{A}^1}^{\trop,\alpha',\beta'}(d-1,g')
  \end{aligned}
	\end{equation}
  for all $ d,g,\alpha,\beta $ as above with $ d>1 $, where the third and fourth sum are
  taken over all $\alpha'$, $\beta'$ and $g'$ satisfying 
  \begin{align*}
    \alpha' &\leq \alpha; \\
    \beta' &\geq \beta; \\
    I\alpha'+I\beta'&=d-1; \\
    g-g'&=|\beta'-\beta|-1; \\
    d-2 &\geq g-g',
  \end{align*}
  and, additonally, $\beta'_i-\beta_i=0$ for all $i$ even in case of the third sum, resp.\ there exists an even $i$ such that $\beta'_i-\beta_i\neq 0$ for the fourth sum.
\end{theorem}

 Note that the sums in Equation~\eqref{eq:Caporaso_Harris} are finite: for the first two summands, this is true since there are only finitely many $\beta_k>0$, as $\beta$ is restricted by $I\alpha+I\beta=d$. For the third and fourth summand, this is true since the possible $\alpha'$ and $\beta'$ are restricted by $I\alpha'+I\beta'=d-1$, and the genus $g'$ is restriced by the fourth condition.

\begin{proof}
We follow the proof of the tropical Caporaso-Harris formula for the complex count, see Theorem 4.3 \cite{GM053}. As the combinatorics of the plane tropical curve count stays the same, we only have to take the necessary changes due to using the arithmetic multiplicity into account.

Let $C$ be a plane tropical curve contributing to $ N_{\mathbb{A}^1}^{\trop,\alpha,\beta}(d,g)$. Pick one of the point conditions and move the point to the very far left. There are two possible outcomes of this procedure. 

\textbf{Case 1:} Moving the point condition places the point on a left end and thus produces a tropical curve $C'$ with one more fixed left end compared to $C$. We show that this case explains the first two summands of the Caporaso-Harris formula~\eqref{eq:Caporaso_Harris}. The new curve $C'$ contributes to $N_{\mathbb{A}^1}^{\trop,\alpha+e_k,\beta-e_k}(d,g)$ for some $k$ with $\beta_k>0$. Let us first assume $k$ is odd, and $\frac{1}{I^{\alpha}\cdot } \mult_\CC(C)$ is odd as well.
Note that $\frac{1}{I^{\alpha}\cdot k} \mult_{\mathbb{C}}(C)$ is automatically odd in this case.
Then the curves $C$ and $C'$ contribute to $ N_{\mathbb{A}^1}^{\trop,\alpha,\beta}(d,g)$ and $N_{\AA^1}^{\trop, \alpha+e_k, \beta-e_k}(d,g)$ with multiplicities
\begin{align*}
	\mult_{\AA^1}^{\alpha, \beta}(C) &= \frac{1}{2} \Big(\frac{1}{I^{\alpha}}  \mult_{\mathbb{C}}(C)-1\Big)\cdot \mathbb{H}+ \big\langle (-1)^i I^\beta \big\rangle \qquad \mbox{ and} \\
	\mult_{\AA^1}^{\alpha + e_k, \beta - e_k}(C') &= \frac{1}{2} \Big(\frac{1}{I^{\alpha}\cdot k} \mult_{\mathbb{C}}(C)-1\Big)\cdot \mathbb{H}+ \Big\langle (-1)^i I^\beta\frac{1}{k} \Big\rangle.
\end{align*}
In the Caporaso-Harris formula~\eqref{eq:Caporaso_Harris} we multiply $\mult_{\AA^1}^{\alpha + e_k, \beta - e_k}(C')$ with $\Big( \frac{k-1}{2}\cdot \mathbb{H}+\langle k \rangle\Big)$, which gives an expression of the form $M\cdot \mathbb{H}+ \big\langle (-1)^i I^\beta \big\rangle$. To determine $M$, we can evaluate at the complex numbers, yielding 
$$2M+1 = \Big(\frac{1}{I^{\alpha}\cdot k} \mult_{\mathbb{C}}(C)\Big)\cdot k,$$ 
hence the product in the Caporaso-Harris formula~\eqref{eq:Caporaso_Harris} recovers $\mult_{\AA^1}^{\alpha, \beta}(C)$, as desired.

Next, we assume $k$ is odd but $\frac{1}{I^{\alpha} } \mult_\CC(C)$ is even. In particular $\frac{1}{I^{\alpha}\cdot k} \mult_\CC(C)$ is even as well. We obtain for the multiplicities of $C$ and $C'$ 
\begin{align*}
	\mult_{\AA^1}^{\alpha, \beta}(C) &= \frac{1}{2} \Big(\frac{1}{I^{\alpha}} \mult_{\mathbb{C}}(C) \Big)\cdot \mathbb{H} \qquad \mbox{ and} \\
	\mult_{\AA^1}^{\alpha + e_k, \beta - e_k}(C') &= \frac{1}{2} \Big(\frac{1}{I^{\alpha}\cdot k} \mult_{\mathbb{C}}(C) \Big)\cdot \mathbb{H}.
\end{align*}
In the Caporaso-Harris formula~\eqref{eq:Caporaso_Harris} we multiply $\mult_{\AA^1}^{\alpha + e_k, \beta - e_k}(C')$ with $\Big( \frac{k-1}{2}\cdot \mathbb{H}+\langle k \rangle\Big)$, which gives an integer multiple of $\mathbb{H}$, which we can again determine by inserting the complex numbers: it is $\big(\frac{1}{I^{\alpha}\cdot k} \mult_{\mathbb{C}}(C)\big)\cdot k$, which equals $\big(\frac{1}{I^{\alpha}} \mult_{\mathbb{C}}(C)\big)$ as desired.

Finally, we assume $k$ is even. Then $\frac{1}{I^{\alpha} } \mult_\CC(C)$ is even as well and
$C$ contributes to $ N_{\mathbb{A}^1}^{\trop,\alpha,\beta}(d,g)$ with arithmetic $(\alpha,\beta)
  $-multiplicity $$\mult_{\AA^1}^{\alpha, \beta}(C) = \frac{1}{2}\Big(\frac{1}{I^{\alpha}}  \mult_{\mathbb{C}}(C)\Big)\cdot \mathbb{H}.$$
However, for $C'$ we have to distinguish cases
\begin{equation*}
	\mult_{\AA^1}^{\alpha + e_k, \beta - e_k}(C') =
	\begin{cases}
		\frac{1}{2} \big(\frac{1}{I^{\alpha}\cdot k} \mult_{\mathbb{C}}(C)-1\big)\cdot \mathbb{H}+ \big\langle (-1)^i  I^\beta\frac{1}{k} \big\rangle & \mbox{if $\frac{1}{I^{\alpha}\cdot k} \mult_\CC(C)$ is odd,} \\
		\frac{1}{2} \big(\frac{1}{I^{\alpha}\cdot k} \mult_{\mathbb{C}}(C)\big)\cdot \mathbb{H} & \mbox{else.}		
	\end{cases}
\end{equation*}
%If we view the end of weight $k$ to be fixed, the curve obtains the arithmetic $(\alpha,\beta)
%  $-multiplicity 
%  $$\frac{1}{2}\cdot \big(\frac{1}{I^{\alpha}\cdot k} \cdot \mult_{\mathbb{C}}(C)\big)\cdot \mathbb{H}, \mbox{ or}$$ 
%   $$\frac{1}{2}\cdot \big(\frac{1}{I^{\alpha}\cdot k} \cdot \mult_{\mathbb{C}}(C)-1\big)\cdot \mathbb{H}+\langle (-1)^i\cdot I^\beta\frac{1}{k}\rangle ,$$
%  depending on whether $\frac{1}{I^{\alpha}\cdot k}\cdot \mult_\CC(C)$ is still even or odd.
  In any case, since we multiply with $\frac{k}{2}\cdot \mathbb{H}$, we obtain an integer multiple of $\mathbb{H}$ which we can determine by inserting the complex numbers. We obtain $\big(\frac{1}{I^{\alpha}\cdot k} \mult_{\mathbb{C}}(C)\big)\cdot k$ in both cases, as desired.
  This completes the first case and explains the first two summands of Formula~\eqref{eq:Caporaso_Harris}.
  
\textbf{Case 2:}
If the point we move to the very left does not fix an additional end, then a floor (i.e.\ a path from a downwards facing end to a diagonal end) splits off on the left, and the curve looks as sketched in Figure~\ref{fig-CHfloor}. We denote the remaining part by $C'$ (which is a plane tropical curve of degree $d-1$) and we compute the total multiplicity of curves $C$ in terms of the multiplicities of all such $C'$. This will explain the remaining two summands of Formula~\eqref{eq:Caporaso_Harris}.

\begin{figure}
	\begin{center}
		\input{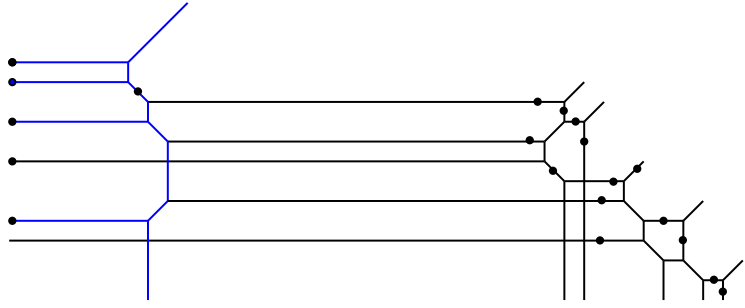}
	\end{center}
	\caption{Moving a point to the very left, a floor (shown in blue) can split off a tropical curve $C$. We denote the rest of the curve (in black) by $C'$.}
	\label{fig-CHfloor}
\end{figure}  

 Assume that $\alpha' \leq
    \alpha$ of the fixed horizontal ends of $C$ only bypass the floor and are not adjacent to a $3$-valent vertex of the floor. Then
    $C'$ has $\alpha'$ fixed horizontal ends. Given a curve $C'$ of degree
    $d-1$ with $\alpha'$ fixed ends which passes through our remaining points,
    there are $\binom{\alpha}{\alpha'}$ possibilities to choose which fixed
    ends of $C$ belong to $C'$. The remaining horizontal ends of $C'$ are non-fixed. Let $\beta'$ be a sequence which fulfills $I\beta'=d-1-
    I\alpha'$, hence a possible choice of weights for the non-fixed ends of
    $C'$. Assume that $\beta''\leq \beta'$ of these ends are adjacent to 
    $3$-valent vertices of the floor whereas $\beta'-\beta''$ ends
    bypass $C\setminus C'$. 
    We claim that every nonfixed end of $C$ is a nonfixed end of $C'$. Indeed, by Lemma~4.20 \cite{Mi03} the one point we moved cannot separate more than two (nonfixed) ends. However, the floor already contains one nonfixed end on either side of the marked point (the down end and the diagonal end). 
    Therefore, all ends of the floor of direction $(-1, 0)$ are fixed.
    So all the $\beta$ nonfixed
    ends of direction $(-1,0)$ have to bypass the floor $C\setminus C'$ --- therefore
    they have to be ends of $C'$. That is, $\beta'-\beta''= \beta$ (in
    particular $\beta' \geq \beta$). Given $C'$, there are $ \binom{\beta'}{
    \beta}$ possibilities to choose which ends of $C'$ are also ends of $C$.
   
Furthermore, note that the complex multiplicity of the floor $C \setminus C'$ is
\[ \mult_\CC(C \setminus C') = I^{\alpha - \alpha'}\cdot I^{\beta' - \beta}. \]
This follows since every triangle dual to a vertex in the floor lies between two lattice columns and accordingly has an area which is equal to the weight of the horizontal adjacent edge.
It follows that the arithmetic $(\alpha-\alpha',\beta'-\beta)$-multiplicity of the floor equals the prefactors in the third resp.\ fourth sum of Formula~\eqref{eq:Caporaso_Harris} which are not combinatorial.  
We have to distinguish the two cases depending on whether $ I^{\beta'-\beta}$ is even or not, as the arithmetic multiplicity of the floor differs for these two cases. 
%This leads to the two summands, the third and the fourth, in our formula.

By Lemma \ref{lem-vertexmult} it follows that the arithmetic $(\alpha,\beta)$-multiplicity of $C$ equals the arithmetic $(\alpha',\beta')$-multiplicity of $C'$ times this prefactor. 

 To determine the genus $g'$ of $C'$, note that $C'$ has by $
    |\alpha+\beta''|$ fewer vertices than $C$ and by $|\alpha+\beta''|-1+
    |\beta''|$ fewer bounded edges --- there are $|\alpha+\beta''|-1$ bounded
    edges in $C\setminus C'$, and $ |\beta''|$ bounded edges are cut. Hence
    $g'=g-(|\beta''|-1)$. Furthermore, $g-g' \leq d-2$ as at most $d-2$ loops may
    be cut. Now given a curve $C'$ with $\alpha'$ fixed and $\beta'$ nonfixed
    bounded edges through the remaining points, and having chosen
    which of the $\alpha$ fixed ends of $C$ are also fixed ends of $C'$ and
    which of the $\beta'$ ends of $C'$ are also ends of $C$, there is only one
    possibility to add a floor through the very left point to make it a
    possible curve $C$ with $\alpha$ fixed ends and $\beta$ nonfixed. The
    positions and directions of all bounded edges are prescribed by the point, by the positions of the $\beta'-\beta$ ends to the left of $C'$, and
    by the $\alpha-\alpha'$ fixed ends. Hence we can count the possibilities
    for $C'$ (times the factor $\binom{\alpha}{\alpha'}\cdot \binom{\beta'}{
    \beta}\cdot I^{\beta'-\beta}$) instead of the possibilities for $C$, where
    the possible choices for $\alpha'$, $\beta'$ and $g'$ have to satisfy 
    the listed conditions.

\end{proof}

%%%%%%%%%%%%%%%%%%%%%%%%%%%%%%%%%%%%%%%%%%%%%%%%%%%%%%%%%%%%%%%%%%%%%%%%%%%%%%%%%%%%%%%%%
\section{Arithmetic multiplicities of floor diagrams and polynomiality properties}
%%%%%%%%%%%%%%%%%%%%%%%%%%%%%%%%%%%%%%%%%%%%%%%%%%%%%%%%%%%%%%%%%%%%%%%%%%%%%%%%%%%%%%%%%

\begin{figure}	
	\begin{center}
		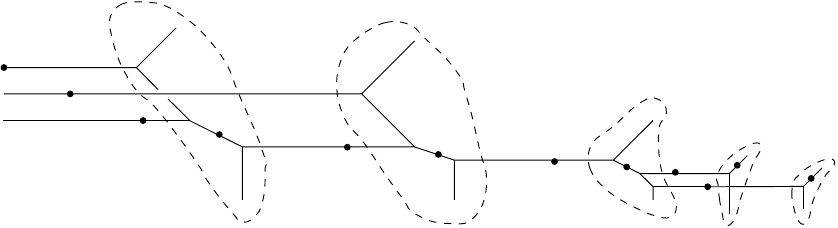
		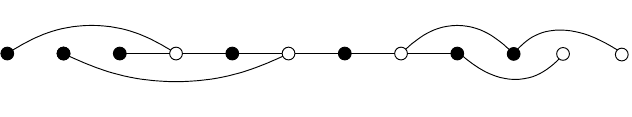
	\end{center}
	\caption{A floor decomposed plane tropical curve of degree $5$ and genus $0$ and its associated floor diagram.}
	\label{fig:floordecomposed}
\end{figure}

The key idea for the proof of the Caporaso-Harris formula in Theorem~\ref{thm-caporaso-harris} was to move one point very far away from the others to the left. This allowed us to either fix a horizontal end with this point or split off a floor from the tropical curve. If we iteratively repeat this process with the remaining points, we can achieve a configuration where all the point conditions lie in a small horizontal strip and thereby obtain a \emph{floor-decomposed} tropical curve. Floor diagrams are a combinatorial gadget encoding the floor-and-elevator structure of a floor-decomposed tropical curve \cite{BM08, FM09}. Given a floor-decomposed tropical curve, its floor diagram is given by shrinking each floor to one point (which we mark as a white vertex to distinguish from the points on horizontal edges, which we mark as black vertices). See Figure \ref{fig:floordecomposed} for an example. 
 There is a bijection between tropical curves through a configuration of points in a small horizontal strip and marked floor diagrams.

We highlight two applications of floor diagrams. Firstly, they were used in the study of piecewise polynomiality properties for counts of curves on Hirzebruch surfaces \cite{AB14}. Secondly, using floor diagrams it was shown \cite{FM09, Blo11} that the number of curves in $\mathbb{P}^2$ with a fixed number of nodal singularities $\delta$ and degree $d$ through an appropriate number of points is eventually polynomial in $d$.%\footnote{In fact, this can be generalized to counts on toric surfaces with \emph{h-transverse} polygon \cite{AB13}, see Definition~\ref{def-hTransverse}.}

In this section we will use the $\mathbb{A}^1$-multiplicities introduced in Definition~\ref{def-troparith2} to define $\mathbb{A}^1$-multi\-plicities for floor diagrams in Definition~\ref{def-arithmultfloor}. 
The $\mathbb{A}^1$-multi\-plicity of a (marked) floor diagram is designed to be equal to the $\mathbb{A}^1$-multi\-plicity of the corresponding tropical curve.

We then revisit the counts of curves in Hirzebruch surfaces and node polynomials to show that the adapted floor diagram count with arithmetic multiplicities satisfies certain polynomiality properties as well. Again the combinatorics of a floor diagram count stays the same no matter with which multiplicity we count.

\subsection{General definition}

\begin{definition} \label{def-hTransverse}
	A lattice polygon $\Delta$ is \emph{h-transverse} if and only if for every $a \in \ZZ$ the intersection 
	\[ \Delta \cap \big\{ (x, y) \in \RR^2 \mathrel{\big |} a \leq x \leq a+1 \big\} \]
	is a (possibly empty) lattice polygon.
\end{definition}

\begin{figure}[h!]
	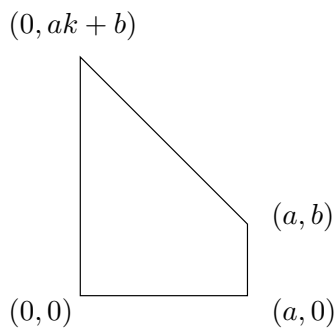
	\caption{The polygon $\Delta_k(a,b)$ defining the Hirzebruch surface $\mathbb{F}_k$ as a toric surface embedded in projective space with hyperplane section the class of a curve of bidegree $(a,b)$. The vertical sides corresponds to the sections $B$ (left) and $E$ (right), see Subsection~\ref{subsec-Hirzebruch}.}
	\label{fig-Hirz}
\end{figure}

Floor diagrams can be defined for any h-transverse lattice polygon \cite{AB13}. 
However, to limit notational complexity we will restrict attention to $\Delta$ of the form in Figure~\ref{fig-Hirz}. Note that $\Delta_d = \conv\big\{ (0,0), (0, d), (d, 0) \big\}$ is a special case of this.

\begin{definition}[Floor diagram]\label{def-floor}
	%Fix some genus $g \geq 0$ and a lattice polygon $\Delta = \Delta_k(a,b)$ as in Figure~\ref{fig-Hirz}.
	Fix an integer $k \geq 0$. 
	Let $\Gamma$ be a graph without loops on the vertex set $\{1, \ldots, a\}$. Then $\Gamma$ is naturally an oriented graph by orienting each edge towards the numerically lager vertex. If $w : E(\Gamma) \to \ZZ_{\geq 1}$ is an edge-weight function, then the pair $(\Gamma, w)$ is a \emph{$\Delta$-floor diagram} if and only if at every vertex $v \in V(\Gamma)$ the \emph{divergence} is
	\[ \operatorname{div}(v) := \sum_{\text{incoming edges }e} w(e) - \sum_{\text{outgoing edges } e} w(e) \leq k. \]
	A floor diagram $\cD$ is said to be \emph{connected} or to have \emph{genus $g$} if and only if this holds for the underlying graph $\Gamma$.
	An edge $i \to i+1$ is called \emph{short} if it has weight 1.
\end{definition}

We fix a polygon as in Figure \ref{fig-Hirz}, a partition $\bw^{(r)}$ of $b$ (the partition of weights of right horizontal ends), and a partition $\bw^{(l)}$ of $ak+b$ (the partition of weights of left horizontal ends).
We require ends of primitive direction $(0,-1)$ and $(1,k)$ to be of weight 1. With this data we fix a Newton fan i.e.\ a degree of tropical curves which we denote henceforth by $\Delta_k(\bw^{(r)},\bw^{(l)},a)$. We also fix a genus $g$.

\begin{definition}[Marking]\label{def-marking}
	Let $\cD = (\Gamma, w)$ be a connected $\Delta$-floor diagram of genus $g$. 
	A \emph{marking} $\cD'$ of $\cD$ of degree $\Delta_k(\bw^{(r)},\bw^{(l)}, a)$ is a connected floor diagram that arises from $\cD$ by the following procedure. We call the vertices of $\cD$ \emph{white} and insert new \emph{black} vertices as follows.
	\begin{enumerate}
		\item Subdivide each edge in $\cD$ by inserting a black vertex in the middle. The new edges inherit weight and orientation.
		\item For each $i = 1, \ldots, \ell(\bw^{(l)})$ add a black vertex $v_i$ and an edge of weight $\bw_i^{(l)}$ starting at $v_i$ and ending in some white vertex. Similarly, for each $j = 1, \ldots, \ell(\bw^{(r)})$ add a black vertex $w_j$ and an edge of weight $\bw_j^{(r)}$ starting at some white vertex and ending in $w_j$. This has to be done such that in the end the divergence of every white vertex is equal to $k$.
		\item Extend the total order on the white vertices to a total order on the set of vertices of $\cD'$ in a way that is compatible with the orientation of the edges.
	\end{enumerate}
	Given a floor diagram $\cD$ with two markings $\cD'$ and $\cD''$, we say that $\cD'$ and $\cD''$ are \emph{equivalent} if there is an isomorphism between the underlying graphs of $\cD'$ and $\cD''$ which respects the weight functions and which fixes the white vertices.
\end{definition}

For a floor diagram $\cD$ we denote the number of markings up to equivalence by $\nu(\cD)$. Typically, a floor diagram count is given as a sum of $\nu(\cD)$, weighted by some multiplicity $\mult(\cD)$, where the sum runs over all floor diagrams with degree and genus as prescribed by the problem at hand.

\begin{definition}[Arithmetic multiplicity of a floor diagram]\label{def-arithmultfloor}
	Let $e$ be an edge of a floor diagram with weight $w$. 
	We define the \emph{arithmetic multiplicity of $e$} to be 
	$$\mult_{\mathbb{A}^1}(e)=
	\begin{cases} 
		\frac{w-1}{2}\cdot \mathbb{H}+\langle w \rangle & \mbox{if $w$ is odd,}\\
		\frac{w}{2}\cdot \mathbb{H} & \mbox{if $w$ is even}.
	\end{cases}$$
	We define the \emph{arithmetic multiplicity} of a floor diagram $\cD$ to be
	$$\mult_{\mathbb{A}^1}(\cD)=\prod_e \mult_{\mathbb{A}^1}(e),$$
	where the product is taken over all edges of $\cD$.
\end{definition}

\subsection{Count of floor diagrams}

Fix a degree $\Delta = \Delta_k(\bw^{(r)},\bw^{(l)},a)$ for tropical curves with Newton polygon of the form depicted in Figure~\ref{fig-Hirz}.

\begin{definition}[Arithmetic count of floor diagrams]\label{def-arithfloorcount}
	Define 
	\[N_{\mathbb{A}^1}^{\floor} \big(\Delta_k(\bw^{(r)},\bw^{(l)},a) , g \big) := \sum_\cD \mult_{\mathbb{A}^1}(\cD) \nu(\cD) \]
	with the sum being taken over all floor diagrams of genus $g$ and $\nu(\cD)$ the number of markings of degree $\Delta_k(\bw^{(r)},\bw^{(l)},a)$ up to equivalence.
\end{definition}

\begin{theorem}\label{thm-floor=trop}
	Fix a genus $g$ and degree $\Delta_k(\bw^{(r)},\bw^{(l)},a)$ for tropical curves with Newton polygon of the form $\Delta_k(a,b)$. 
  Let $K$ be a field of characteristic larger than the diameter of $\Delta_k(\bw^{(r)},\bw^{(l)},a)$.
 Then the arithmetic count of floor diagrams equals the arithmetic count of tropical curves,
	$$N_{\mathbb{A}^1}^{\floor} \big(\Delta_k(\bw^{(r)},\bw^{(l)},a),g \big) = N_{\mathbb{A}^1}^{\trop} \big(\Delta_k(\bw^{(r)},\bw^{(l)},a),g \big)\in \GW(K).$$
\end{theorem}
\begin{proof}
	As proved e.g.\ in Theorem 3.7 \cite{FM09}, there is a bijection between floor diagrams and tropical curves of the respective degree and genus. We only need to check that the multiplicity of the floor diagram equals the multiplicity of the corresponding tropical curve. Let $v$ be a vertex of such a tropical curve and $\delta$ the corresponding triangle in the dual Newton subdivision. Since the tropical curve is floor decomposed, $\delta$ lies between two lattice columns and accordingly the area of $\delta$ is equal to the weight of the horizontal edge $e$ adjacent to $v$. Note that $e$ is an edge in the floor diagram as well. Comparing Definitions~\ref{def-troparith2} and~\ref{def-arithfloorcount} we see that
	\[ \mult_{\AA^1}(v) = \mult_{\AA^1}(e) . \]
	An end $e$ of the tropical curve of weight $w$ has only one incident vertex $v$, which then contributes a factor of $\mult_{\AA^1}(v)$. This end $e$ corresponds to an edge in the marked floor diagram which we denote by $e$ as well and which was added to the floor diagram in the marking process, step (2). This edge contributes a factor of $\mult_{\AA^1}(e) = \mult_{\AA^1}(v)$ to the floor diagram multiplicity as well. Each bounded horizontal edge $e$ of the tropical curve of weight $w$ has two incident vertices $v_1$, $v_2$, contributing a factor of $\mult_{\AA^1}(v_1) = \mult_{\AA^1}(v_2) = \mult_{\AA^1}(e)$ each. On the floor diagram side, $e$ corresponds to an edge that gets subdivided by a black vertex in the marking process, step (1). The two parts of the subdivided edge each contribute $\mult_{\AA^1}(e)$ again.
	From this the claim follows.
	%		
	%This follows as in the proof of the Caporaso-Harris formula since every triangle dual to a vertex in the floor lies between two lattice columns of the dual Newton subdivision and accordingly has an area which is equal to the weight of the horizontal adjacent edge. 
	%We have to distinguish the two cases depending on whether this weight is even or not, as the arithmetic multiplicity of the vertex differs for these two cases. By Lemma \ref{lem-vertexmult}, the arithmetic multiplicity of the whole curve equals the product of its arithmetic vertex multiplicities.
\end{proof}

\subsection{Hirzebruch surfaces}
\label{subsec-Hirzebruch}

For $k \geq 0$, the \emph{Hirzebruch surface} $\mathbb{F}_k$ is defined to be the surface $\mathbb{P}(\mathcal{O}_{\mathbb{P}^1}\oplus \mathcal{O}_{\mathbb{P}^1}(k))$; 
it is a smooth projective toric surface.
The Picard group of $\mathbb{F}_k$ is isomorphic to $\ZZ\times \ZZ$ generated by the classes of the zero section $B$ and a fiber $F$. In particular, the infinity section $E$ satisfies $E = B-kF$. A curve in $\mathbb{F}_k$ has \emph{bidegree $(a,b)$} if its class is $aB+bF$. The polygon depicted in Figure \ref{fig-Hirz} defines $\mathbb{F}_k$ as a projective toric surface polarized by an $(a,b)$ curve. Double Gromov-Witten invariants of Hirzebruch surfaces were introduced in \cite{AB14}. These are counts of curves in a Hirzebruch surface which satisfy, in addition to the right number of point conditions, tangency conditions to the zero and infinity section. The partition of tangency orders along each section is recorded in a tuple of partitions which we denote by $(\bw^{(r)},\bw^{(l)})$. Using the correspondence theorem \cite[Theorem C]{Ran15}, the count of such curves equals the count of tropical curves of degree dual to the polygon in Figure \ref{fig-Hirz}, where the weights of the left and right ends are fixed by  the two partitions $(\bw^{(r)},\bw^{(l)})$. Using Theorem \ref{thm-floor=trop}, this number can (also when enriched via arithmetic multiplicities) be determined using floor diagrams.

Let us fix lengths of partitions $\ell_r, \ell_l \geq 1$. 
We denote by $H$ the set of partitions $(\bw^{(r)},\bw^{(l)})\in \mathbb{N}^{\ell_r+\ell_l}$ such that $\sum \bw^{(l)}_i = ak+b$ and $\sum \bw^{(r)}_j = b$. We also require that all entries of $\bw^{(r)}$ and $\bw^{(l)}$ are odd, to avoid cases distinctions.
Then the arithmetic count of floor diagrams can be viewed as a function
\begin{equation} \label{eq:floor_diag_count}
    \begin{aligned}
    	N_{\mathbb{A}^1}^{\floor}(k,a,g): \qquad H\qquad &\longrightarrow \GW(K) \\ 
    	(\bw^{(r)},\bw^{(l)}) &\longmapsto N_{\mathbb{A}^1}^{\floor}\big(\Delta_k(\bw^{(r)},\bw^{(l)},a),g \big).
    \end{aligned}
\end{equation}
We can see that the arithmetic count of floor diagrams (and with it, using Theorem \ref{thm-floor=trop} the arithmetic count of tropical curves of degree $\Delta_k(\bw^{(r)},\bw^{(l)},a)$ and genus $g$ is piecewise polynomial in the following sense:

\begin{theorem} \label{thm-piecepoly}
	Let $K$ be a field of characteristic $0$. The function $N_{\mathbb{A}^1}^{\floor}(k,a,g)$ associates to a tuple $(\bw^{(r)},\bw^{(l)})$ of odd weights an element in $\GW(K)$ which equals a multiple of $\mathbb{H}$ plus a sum of quadratic forms of the form $\big\langle (-1)^i\cdot \prod\bw^{(r)}_i\prod\bw^{(l)}_j \big\rangle$. The coefficient of  $\mathbb{H}$ is piecewise polynomial in the entries of $\bw^{(r)}$ and $\bw^{(l)}$.
\end{theorem}
\begin{proof}
	For an individual floor diagram, corresponding to the tropical curve $C$, we obtain an arithmetic multiplicity which is of the form
	$\frac{\mult_\CC(C)-1}{2}\cdot \mathbb{H}+\big\langle (-1)^i\cdot \prod\bw^{(r)}_i\prod\bw^{(l)}_j \big\rangle$ resp.\
	$\frac{\mult_\CC(C)}{2}\cdot \mathbb{H}$, depending on whether it contains even edges or not.
	The statement then follows since the complex multiplicity is piecewise polynomial by Theorem 1.3 \cite{AB14}. The piecewise structure appears because not all floor diagrams exist for every choice of $(\bw^{(r)},\bw^{(l)})$. 
\end{proof}
The idea to see the polynomiality of the complex multiplicity can be described quickly and we include it for the sake of completeness: in the case of genus $0$, the divergence condition fixes the weights of all edges of a floor diagram to be affine-linear forms in the weights of the ends $\bw^{(r)}_i,\bw^{(l)}_j$. Their product is, accordingly, a polynomial in the $\bw^{(r)}_i$ and $\bw^{(l)}_j$. In the case of higher genus, we can fix a set of $g$ edges separating all cycles and use variables $i_1,\ldots,i_g$ for their weights. The weight of any other edge is then an affine-linear form in the $\bw^{(r)}_a,\bw^{(l)}_b$ and the new variables $i_c$. Summing over all possibilities for the $i_c$ amounts to computing Bernoulli formulas which are also polynomial in the $\bw^{(r)}_a,\bw^{(l)}_b$.

\subsection{Node polynomials}

Let us now discuss the \emph{Severi degree} $N^{d, \delta}$ for some degree $d \geq 1$ and some \emph{cogenus} (i.e. number of nodal singularities) $\delta \geq 0$. This quantity $N^{d, \delta}$ is defined as the number of (possibly reducible) curves in $\mathbb{P}^2$ of degree $d$ with $\delta$ many nodal singularities, which pass through $\frac{d(d+3)}{2}-\delta$ points in general position. 
%Since in the plane $g = \frac{(d-1)(d-2)}{2} - \delta$, the Severi degree can also be understood as a Gromov-Witten invariant. 

Recall that plane curves of degree $d$ have Newton polygon $\Delta_d = \conv\big\{  (0,0), (0, d), (d, 0)\big\}$, which is a special case of Figure~\ref{fig-Hirz}. Therefore, the notions of (marked) floor diagrams apply to this case as well. We specialize even further by taking the partition of weights of left ends to be $d =1 + \cdots + 1$ (obviously, there are no right ends in $\Delta_d$).

\begin{definition}
	For fixed $d$ and $\delta$ we define the\emph{ arithmetic count of floor diagrams }to be
	\begin{equation} \label{eq-floorDiagramCount}
		N_{\mathbb{A}^1}^{\delta}(d) := \sum_{\cD} \mult_{\mathbb{A}^1}(\cD) \nu(\cD)\in \GW(K),
	\end{equation}
	where the sum runs over all floor diagrams of degree $d$ and cogenus $\delta(\cD) := \frac{(d-1)(d-2)}{2} - g(\cD)$ equal to $\delta$. 
\end{definition}

 Using Theorem \ref{thm-floor=trop}, $N_{\mathbb{A}^1}^{\delta}(d)$ equals the arithmetic count of tropical curves of degree $d$ and cogenus $\delta$ passing through the right number of points. Using the techniques of the Correspondence Theorem \ref{thm-corres-new}, $N_{\mathbb{A}^1}^{\delta}(d)$ equals the corresponding count of algebraic curves close to the tropical limit.
 
As an example, adapting the computations from \cite{BG14}, we have that if $d\geq\delta$, then
\begin{align*}
    N_{\mathbb{A}^1}^0(d) &= \langle 1 \rangle,\\
    N_{\mathbb{A}^1}^1(d) &= (d-1)(d-2)\mathbb{H}+
    \left(d^2-1\right) \langle 1 \rangle,\\
    N_{\mathbb{A}^1}^2(d) &= \left(2d^4-9d^3+4d^2+21d-18\right)\mathbb{H}+
    \left(\frac{d^4-4d^2-3d+6}{2}\right)\langle 1 \rangle.
\end{align*}
The main result of this section is the following.

\begin{theorem}\label{thm-polynomiality}
	Let $K$ be a field of characteristic $0$. For every $\delta$ there are polynomials $P, Q \in \mathbb{Q}[d]$ of degree $2\delta$ such that for all sufficiently large $d$ we have 
	\[ N_{\mathbb{A}^1}^\delta(d) = P(d) \mathbb{H} + Q(d) \langle 1 \rangle. \]
\end{theorem}

Fomin and Mikhalkin \cite{FM09} have decomposed marked floor diagrams into elementary building blocks, called \emph{templates}. This decomposition is key for proving Theorem~\ref{thm-polynomiality} and for the readers' convenience we repeat the definition here.

\begin{definition}[Templates]
	A \emph{template} is a (possibly disconnected) loop-less graph $\Gamma$, defined on a linearly ordered vertex set $\{0, \ldots, l\}$, together with an edge-weight function $w : E(\Gamma) \to \ZZ_{\geq 1}$ such that
	\begin{enumerate}
		\item the weight of an edge from any vertex $i$ to $i+1$ is not 1 (i.e. there are no short edges) and
		\item for every $1\leq i<l$ there is an edge $j \to k$ with $j < i < k$.
	\end{enumerate} 
	The \emph{length} of a template $\Gamma$ is $\ell(\Gamma) =l$.
\end{definition}

Every marked floor diagram $\cD$ can be decomposed into a sequence of templates as follows. First, construct the \emph{modification} $\cD'$ of $\cD$ by adding one more white vertex named $0$ at the very left of $\cD$ and for every vertex $i$ of $\cD$ add edges $0 \to i$ of weight 1 until all vertices other than $0$ have divergence equal to 1.
Next, we decompose $\cD'$ into templates by deleting all edges of the form $i \to i+1$ and weight 1. What remains is an ordered sequence of templates $\Gamma_1, \ldots, \Gamma_m$ where two neighboring $\Gamma_j$ and $\Gamma_{j+1}$ overlap in at most one vertex. If we additionally remember the index $k_j$ of the vertex in $\cD'$ where the template $\Gamma_j$ starts, then this process can be reversed. Indeed, we simply need to align the templates according to their start indices, re-insert the vertices of $\cD'$ that are not part of any template, and then add short edges $i \to i+1$ at every vertex $i$ until the divergence is equal to 1 everywhere. This recovers the modification $\cD'$ and in order to recover $\cD$ we delete the left-most vertex $0$ from $\cD'$ and all edges starting there.

This reasoning suggests that in the floor diagram count \eqref{eq-floorDiagramCount} we could just as well take the sum over sequences of templates $\Gamma_1, \ldots, \Gamma_m$ together with \emph{valid} sequences of start indices $k_1, \ldots, k_m$. A sequence of start indices is valid for this purpose if the following hold.
\begin{enumerate}
	\item The templates do not overlap in more than one vertex, i.e. $k_{j+1} \geq k_j + \ell(\Gamma_j)$ for all $j = 1, \ldots, m-1$.
	\item The last template $\Gamma_m$ in the sequence fits at the end of the floor diagram $\cD'$, i.e. $k_m + \ell(\Gamma_m) \leq d$.
	\item By construction, the vertex 0 in $\cD'$ has only outgoing edges of weight 1. This needs to be respected by $\Gamma_1$, i.e. we need to have $k_1 \geq 1$ if the first vertex of $\Gamma_1$ has an outgoing edge of weight $>1$.
	\item In $\cD'$ the vertex $i$ has precisely $d-i$ outgoing edges when counted with their weights. Hence, if a vertex in a template $\Gamma_j$ has outgoing edges of total weight $d-i$, then the start index $k_j$ of $\Gamma_j$ has to be chosen small enough for that vertex to land at most at position $i$ in $\cD'$. Here, edges bypassing a vertex also have to be counted as \enquote{outgoing} since
	\tikzset{every picture/.style={line width=0.75pt}} %set default line width to 0.75pt        
	\begin{tikzpicture}[x=0.75pt,y=0.75pt,yscale=-1,xscale=1]
		%uncomment if require: \path (0,784); %set diagram left start at 0, and has height of 784
		
		%Shape: Circle [id:dp9036089823928418] 
		\draw  [fill={rgb, 255:red, 0; green, 0; blue, 0 }  ,fill opacity=1 ] (138.66,420.11) .. controls (138.66,419.28) and (139.33,418.61) .. (140.16,418.61) .. controls (140.99,418.61) and (141.66,419.28) .. (141.66,420.11) .. controls (141.66,420.94) and (140.99,421.61) .. (140.16,421.61) .. controls (139.33,421.61) and (138.66,420.94) .. (138.66,420.11) -- cycle ;
		%Shape: Circle [id:dp37664973843101957] 
		\draw  [fill={rgb, 255:red, 0; green, 0; blue, 0 }  ,fill opacity=1 ] (158.71,420.21) .. controls (158.71,419.39) and (159.39,418.71) .. (160.21,418.71) .. controls (161.04,418.71) and (161.71,419.39) .. (161.71,420.21) .. controls (161.71,421.04) and (161.04,421.71) .. (160.21,421.71) .. controls (159.39,421.71) and (158.71,421.04) .. (158.71,420.21) -- cycle ;
		%Shape: Circle [id:dp5373187623295002] 
		\draw  [fill={rgb, 255:red, 0; green, 0; blue, 0 }  ,fill opacity=1 ] (178.57,420.21) .. controls (178.57,419.39) and (179.24,418.71) .. (180.07,418.71) .. controls (180.9,418.71) and (181.57,419.39) .. (181.57,420.21) .. controls (181.57,421.04) and (180.9,421.71) .. (180.07,421.71) .. controls (179.24,421.71) and (178.57,421.04) .. (178.57,420.21) -- cycle ;
		%Curve Lines [id:da9761213194967668] 
		\draw    (140.16,420.11) .. controls (150.14,410.68) and (170.57,410.54) .. (180.07,420.21) ;
	\end{tikzpicture} 
	fits at some index in $\cD'$ if and only if 
	\tikzset{every picture/.style={line width=0.75pt}} %set default line width to 0.75pt        
	\begin{tikzpicture}[x=0.75pt,y=0.75pt,yscale=-1,xscale=1]
		%uncomment if require: \path (0,784); %set diagram left start at 0, and has height of 784
		
		%Shape: Circle [id:dp9036089823928418] 
		\draw  [fill={rgb, 255:red, 0; green, 0; blue, 0 }  ,fill opacity=1 ] (138.66,420.11) .. controls (138.66,419.28) and (139.33,418.61) .. (140.16,418.61) .. controls (140.99,418.61) and (141.66,419.28) .. (141.66,420.11) .. controls (141.66,420.94) and (140.99,421.61) .. (140.16,421.61) .. controls (139.33,421.61) and (138.66,420.94) .. (138.66,420.11) -- cycle ;
		%Shape: Circle [id:dp37664973843101957] 
		\draw  [fill={rgb, 255:red, 0; green, 0; blue, 0 }  ,fill opacity=1 ] (158.71,420.21) .. controls (158.71,419.39) and (159.39,418.71) .. (160.21,418.71) .. controls (161.04,418.71) and (161.71,419.39) .. (161.71,420.21) .. controls (161.71,421.04) and (161.04,421.71) .. (160.21,421.71) .. controls (159.39,421.71) and (158.71,421.04) .. (158.71,420.21) -- cycle ;
		%Shape: Circle [id:dp5373187623295002] 
		\draw  [fill={rgb, 255:red, 0; green, 0; blue, 0 }  ,fill opacity=1 ] (178.57,420.21) .. controls (178.57,419.39) and (179.24,418.71) .. (180.07,418.71) .. controls (180.9,418.71) and (181.57,419.39) .. (181.57,420.21) .. controls (181.57,421.04) and (180.9,421.71) .. (180.07,421.71) .. controls (179.24,421.71) and (178.57,421.04) .. (178.57,420.21) -- cycle ;
		%Curve Lines [id:da9761213194967668] 
		\draw    (140.16,420.11) .. controls (146,415.54) and (154.43,415.4) .. (160.21,420.21) ;
		%Curve Lines [id:da21118773198649465] 
		\draw    (160.01,420.11) .. controls (165.86,415.54) and (174.29,415.4) .. (180.07,420.21) ;
	\end{tikzpicture}
	fits at that index.
	As a consequence, one can define a maximal start index $k_{\textrm{max}}(\Gamma_j)$ for each template $\Gamma_j$ and then we need to have $k_j \leq k_{\textrm{max}}(\Gamma_j)$.
\end{enumerate}
With this characterization of valid start indices, we can rewrite the sum \eqref{eq-floorDiagramCount} in terms of sequences of templates. Let us now describe how multiplicity, cogenus, and number of markings of a floor diagram $\cD$ break down into contributions from the templates. 

For a template $\Gamma$ define its \emph{arithmetic multiplicity}
\[ \mult_{\mathbb{A}^1} (\Gamma) := \prod_{e} \mult_{\mathbb{A}^1} (e) \]
with the product being taken over all edges in $\Gamma$. If now $\cD$ decomposes into a sequence of templates $\Gamma_1, \ldots, \Gamma_m$, then $\mult_{\mathbb{A}^1}(\cD)$ is given as the product of the multiplicities $\mult_{\mathbb{A}^1}(\Gamma_j)$.

To see how the cogenus $\delta(\cD)$ decomposes, consider the reference floor diagram $\cE$ with cogenus $\delta(\cE) = 0$ consisting only of short edges. Then any floor diagram $\cD$ can be obtained from $\cE$ by successively replacing sets of $w$ many sequences of short edges $i \to i+1 \to \cdots \to j$ with a single edge $i \to j$ of weight $w$. Each such step increases the cogenus by $(j-i)w - 1$, i.e. the cogenus of $\cD$ can be expressed as 
$\delta(\cD) = \sum_{e = i \to j} \big((j - i)w(e) - 1 \big)$.
If we define the cogenus of a template $\Gamma$ to be 
\[\delta(\Gamma) := \sum_{e = i \to j} \big( (j-i)w(e) - 1 \big),\] 
then it is clear that $\delta(\cD)$ is the sum of the cogenera $\delta(\Gamma_j)$ taken over the template decomposition $\Gamma_1, \ldots, \Gamma_m$ of $\cD$.

Concerning the number of markings it is easy to see that $\nu(\cD) = \nu(\cD')$. Now if we let $\Gamma_1, \ldots, \Gamma_m$ once more be the template decomposition of $\cD$ then for each $\Gamma_j$, the number of markings of the subgraph of $\cD'$ that is spanned by $\Gamma_j$ depends only on $\Gamma_j$ and its start index $k_j$. We denote this quantity by $\nu(\Gamma_j, k_j)$. Note that the subgraphs of $\cD'$ connecting the end vertex of $\Gamma_{j-1}$ and the start vertex of $\Gamma_j$ consist exclusively of parallel short edges, i.e. up to equivalence there is a unique way of marking them. Consequently
\[ \nu(\cD) = \prod_{j = 1}^m \nu(\Gamma_j, k_j) \ . \]
These above observations constitute a proof for the following lemma.

\begin{lemma}[compare Equation~(5.13)~\cite{FM09}] For all $d$ and $\delta$ we have
	\begin{equation} \label{eq-templateDecomposition}
		N_{\mathbb{A}^1}^\delta(d) = 
		\sum_{\substack{\text{templates } \Gamma_1, \ldots, \Gamma_m \\ \text{such that } \sum \delta(\Gamma_j) = \delta}} 
		\Bigg(\prod_{j = 1}^m\mult_{\mathbb{A}^1}(\Gamma_j) \Bigg)
		\underbrace{
			\sum_{\substack{k_1, \ldots, k_m \\ \text{valid starting indices}}}
			\prod_{j = 1}^m \nu(\Gamma_j, k_j)
		}_{(*)} \ .
	\end{equation}
\end{lemma}

\begin{proof}[Proof of Theorem~\ref{thm-polynomiality}]
	It was shown in the proof of Theorem~5.1 \cite{FM09} that the quantity $(*)$ in Equation~\eqref{eq-templateDecomposition} is a polynomial of degree $\leq 2\delta$. Furthermore, the number of summands in the first sum of Equation~\eqref{eq-templateDecomposition} does not depend on $d$. Finally, we note that $m \leq \delta$ and hence the claim follows from the fact that any product of elements from $\GW(K)$ of the form $a\langle 1\rangle + b \mathbb{H}$ is again of the same form.	
\end{proof}

\bibliographystyle{plain} 
\bibliography{bibliographie}

\end{document}